\documentclass[10pt]{article}
\usepackage{bm}
\usepackage[mathscr]{eucal}
\usepackage{mathrsfs}
\usepackage{amssymb}
\usepackage{calligra}
\usepackage{fontenc}
\usepackage{amsbsy} 
 \usepackage{paralist}
\usepackage{graphicx}
\graphicspath{%
    {converted_graphics/}
    {/}
}
\begin{document}
\newcommand{\ts}{{\tilde{\sf s}}}
\newcommand{\sfv}{{\sf v}}
\newcommand{\sfw}{{\sf w}}
\newcommand{\simge}{\ba{cc}\vspace*{-2.4mm}>\\ \sim\ea }
\newcommand{\simle}{\ba{cc}\vspace*{-2.4mm}<\\ \sim\ea }
\newcommand{\Cdot}{\!\cdot\!}
\newcommand{\sq}{{$\sqcap\!\!\!\!\sqcup$}}
\newcommand{\Eu}{{\rm I\,\!\! E}}
\newcommand{\Io}{\Int{\Omega}{}}
\newcommand{\Id}{\Int{\cald}{}}
\newcommand{\Div}{\mbox{\rm div}\,}
\newcommand{\tr}{\mbox{\rm tr}\,}
\newcommand{\grad}{\mbox{\rm grad}\,}
\newcommand{\supp}{\mbox{\rm supp}\,}
\newcommand{\curl}{\mbox{\rm curl}\,}
\newcommand{\Ido}{\Int{\partial\Omega}{}}
\newcommand{\IdS}{\Int{\Sigma}{}}
\newcommand{\Oint}[2]{{\displaystyle \oint_{#1}^{#2}}}
\newcommand{\Int}[2]{{\displaystyle \int_{ #1}^{ #2}}}
\newcommand{\Lim}[1]{{\displaystyle \lim_{ #1}}}
\newcommand{\Limsup}[1]{{\displaystyle \limsup_{\footnotesize #1}}}
\newcommand{\Liminf}[1]{{\displaystyle \liminf_{\footnotesize #1}}}
\newcommand{\Sup}[1]{{\displaystyle \sup_{#1}}}
\newcommand{\Inf}[1]{{\displaystyle \inf_{#1}}}
\newcommand{\Max}[1]{{\displaystyle \max_{#1}}}
\newcommand{\Min}[1]{{\displaystyle \min_{#1}}}
\newcommand{\Sum}[2]{{\displaystyle \sum_{#1}^{#2}}}
\newcommand{\Prod}[2]{{\displaystyle \prod_{#1}^{#2}}}
\newcommand{\BCup}[2]{{\displaystyle \bigcup_{#1}^{#2}}}
\newcommand{\BCap}[2]{{\displaystyle \bigcap_{#1}^{#2}}}
\newcommand{\Frac}[2]{\displaystyle{\frac{\displaystyle{#1}}{\displaystyle{#2}}}}
\newcommand{\norm}[1]{\left\|{#1}\right\|}
\newcommand{\Norm}[1]{\langle\langle{#1}\rangle\rangle_q}
\newcommand{\No}[1]{\langle\!\langle{#1}\rangle\!\rangle}
\newcommand{\NO}[1]{{\langle{#1}\rangle}_{\lambda,q}}
\newcommand{\beea}{\begin{eqnarray}}
\newcommand{\eeea}{\end{eqnarray}}
\newcommand{\ms}{\medskip\smallskip}
\newcommand{\bs}{\bigskip}
\newcommand{\ps}{\par\smallskip}
\newcommand{\bfe}{{\mbox{\boldmath $e$}} }
\newcommand{\pni}{{\par\noindent}}
\newcommand{\bfq}{{\mbox{\boldmath $q$}} }
\newcommand{\bfz}{{\mbox{\boldmath $z$}} }
\newcommand{\0}{{\mbox{\boldmath $0$}} }
\newcommand{\LE}{\!\!\!&\le&\!\!\!}
\newcommand{\BL}[1]{{\par\smallskip{\bf Lemma #1.}}}
\newcommand{\BT}[1]{{\par\smallskip{\bf Theorem #1.}}}
\newcommand{\Ln}{[\!|}
\newcommand{\Rn}{|\!]}
\newcommand{\n}[1]{{\Ln{#1}\Rn}} 
\newcommand{\nq}[1]{{\Ln{#1}\Rn}_{q}} 
\newcommand{\nqr}[1]{{\Ln{#1}\Rn}_{q,r}} 
\newcommand{\Nq}[1]{{\langle{#1}\rangle}_{q}} 
\newcommand{\Nql}[1]{{\langle{#1}\rangle}_{\lambda,q}} 
\newcommand{\Nqr}[1]{{\langle{#1}\rangle}_{q,r}}
\newcommand{\N}[1]{{|\!\!|\!\!|\,{#1}\,|\!|\!\!|_2}}
\newcommand{\EA}[2]{$$#1$$%
\vspace{-6.mm}
\begin{equation}
\end{equation}
\vspace{-6.mm}
$$
#2
\setlength{\belowdisplayskip}{3mm}
\setlength{\belowdisplayshortskip}{3mm}
$$
}
\newcommand{\A}[2]{$$#1$$%
\vspace{-4.mm}
$$
#2
\setlength{\belowdisplayskip}{3mm}
\setlength{\belowdisplayshortskip}{3mm}
$$
}
\newcommand{\BF}{\begin{footnotesize}}
\newcommand{\EF}{\end{footnotesize}}
\setlength{\jot}{.15in}
\newcommand{\pde}[2]{{\displaystyle \frac{\mbox{$\partial #1$}}{\mbox{$\partial #2$}}}}
\newcommand{\ode}[2]{{\displaystyle \frac{\mbox{$d #1$}}{\mbox{$d #2$}}}}
\newcommand{\f}[2]{\frac{\mbox{$#1$}}{\mbox{$ #2$}}}
\newcommand{\bi}{\begin{itemize}}
\newcommand{\ei}{\end{itemize}}
\newcommand{\ed}{\end{document}}
\newcommand{\be}{\begin{equation}}
\newcommand{\ba}{\begin{array}}
\newcommand{\ea}{\end{array}}
\newcommand{\ee}{\end{equation}}
\newcommand{\eeq}[1]{\label{eq:#1}\end{equation}}
\newcommand{\real}{{\mathbb R}}
\newcommand{\compl}{{\mathbb C}}
\def\Id{\mbox{\boldmath $1$}}
\def\zero{\mbox{\boldmath $0$}}
\newcommand{\PP}{{\rm I\!\!\,P}}
\newcommand{\nat}{{\mathbb N}}
\newcommand{\bfpsi}{\mbox{\boldmath $\psi$}}
\newcommand{\bfchi}{\mbox{\boldmath $\chi$}}
\newcommand{\bfomega}{\mbox{\boldmath $\omega$}}
\newcommand{\bfome}{\mbox{\boldmath $\varpi$}}
\newcommand{\bfvaromega}{\mbox{\boldmath $\varpi$}}
\newcommand{\bfOmega}{\mbox{\boldmath $\Omega$}}
\newcommand{\bfTheta}{\mbox{\boldmath $\Theta$}}
\newcommand{\bfxi}{\mbox{\boldmath $\xi$}}
\newcommand{\bfmu}{\mbox{\boldmath $\mu$}}
\newcommand{\bfx}{\mbox{\boldmath $x$}}
\newcommand{\bfy}{\mbox{\boldmath $y$}}
\newcommand{\bfPsi}{\mbox{\boldmath $\Psi$}}
\newcommand{\bfphi}{\mbox{\boldmath $\varphi$}}
\newcommand{\bfhi}{\mbox{\boldmath $\phi$}}
\newcommand{\bfPhi}{\mbox{\boldmath $\Phi$}}
\newcommand{\bfv}{{\mbox{\boldmath $v$}} }
\newcommand{\bfu}{{\mbox{\boldmath $u$}} }
\newcommand{\bfsf}{{\mbox{\footnotesize\boldmath $s$}} }
\newcommand{\bfuf}{{\mbox{\footnotesize\boldmath $u$}} }
\newcommand{\bfw}{{\mbox{\boldmath $w$}} }
\newcommand{\bff}{{\mbox{\boldmath $f$}} }
\newcommand{\bfa}{{\mbox{\boldmath $a$}} }
\newcommand{\bfi}{{\mbox{\boldmath $i$}} }
\newcommand{\bfj}{{\mbox{\boldmath $j$}} }
\newcommand{\bfc}{{\mbox{\boldmath $c$}} }
\newcommand{\bfo}{{\mbox{\boldmath $o$}} }
\newcommand{\bfp}{{\mbox{\boldmath $p$}} }
\newcommand{\bfkp}{{\mbox{\footnotesize{\boldmath $k$}}} }
\newcommand{\bfka}{{\mbox{\footnotesize{\boldmath $k^*$}}} }
\newcommand{\bft}{{\mbox{\boldmath $t$}} }
\newcommand{\bfd}{{\mbox{\boldmath $d$}} }
\newcommand{\bfl}{{\mbox{\boldmath $l$}} }
\newcommand{\bfr}{{\mbox{\boldmath $r$}} }
\newcommand{\bfk}{{\mbox{\boldmath $k$}} }
\newcommand{\bfA}{{\mbox{\boldmath $A$}} }
\newcommand{\bfS}{{\mbox{\boldmath $S$}} }
\newcommand{\bfO}{{\mbox{\boldmath $O$}} }
\newcommand{\bfM}{{\mbox{\boldmath $M$}} }
\newcommand{\bfP}{{\mbox{\boldmath $P$}} }
\newcommand{\bfB}{{\mbox{\boldmath $B$}} }
\newcommand{\bfR}{{\mbox{\boldmath $R$}} }
\newcommand{\bfC}{{\mbox{\boldmath $C$}} }
\newcommand{\bfD}{{\mbox{\boldmath $D$}} }
\newcommand{\bfQ}{{\mbox{\boldmath $Q$}} }
\newcommand{\bfZ}{{\mbox{\boldmath $Z$}} }
\newcommand{\bfG}{{\mbox{\boldmath $G$}} }
\newcommand{\bfE}{{\mbox{\boldmath $E$}} }
\newcommand{\bfX}{{\mbox{\boldmath $X$}} }
\newcommand{\bfY}{{\mbox{\boldmath $Y$}} }
\newcommand{\bfH}{{\mbox{\boldmath $H$}} }
\newcommand{\bfI}{{\mbox{\boldmath $I$}} }
\newcommand{\bfJ}{{\mbox{\boldmath $J$}} }
\newcommand{\bfN}{{\mbox{\boldmath $N$}} }
\newcommand{\bfh}{{\mbox{\boldmath $h$}} }
\newcommand{\bfm}{{\mbox{\boldmath $m$}} }
\newcommand{\bfone}{{\mbox{\boldmath $1$}} }
\newcommand{\hs}{{\rm I}\!\!\,{\rm R}^3_+}
\newcommand{\cala}{{\cal A}}
\newcommand{\calb}{{\cal B}}
\newcommand{\calc}{{\cal C}}
\newcommand{\cald}{{\cal D}}
\newcommand{\cale}{{\cal E}}
\newcommand{\calf}{{\cal F}}
\newcommand{\calg}{{\cal G}}
\newcommand{\calh}{{\cal H}}
\newcommand{\cali}{{\cal I}}
\newcommand{\calj}{{\cal J}}
\newcommand{\calk}{{\cal K}}
\newcommand{\call}{{\cal L}}
\newcommand{\calm}{{\cal M}}
\newcommand{\caln}{{\cal N}}
\newcommand{\calo}{{\cal O}}
\newcommand{\calp}{{\cal P}}
\newcommand{\calq}{{\cal Q}}
\newcommand{\calr}{{\cal R}}
\newcommand{\cals}{{\cal S}}
\newcommand{\calt}{{\cal T}}
\newcommand{\calu}{{\cal U}}
\newcommand{\calv}{{\cal V}}
\newcommand{\calx}{{\cal X}}
\newcommand{\caly}{{\cal Y}}
\newcommand{\calw}{{\cal W}}
\newcommand{\calz}{{\cal Z}}
\newcommand{\bfsigma}{\mbox{\boldmath $\sigma$}}
\newcommand{\bfSigma}{\mbox{\boldmath $\Sigma$}}
\newcommand{\bftau}{\mbox{\boldmath $\tau$}}
\newcommand{\bfeta}{\mbox{\boldmath $\eta$}}
\newcommand{\bfT}{{\mbox{\boldmath $T$}} }
\newcommand{\bfV}{{\mbox{\boldmath $V$}} }
\newcommand{\bfU}{{\mbox{\boldmath $U$}} }
\newcommand{\bfW}{{\mbox{\boldmath $W$}} }
\newcommand{\bfF}{{\mbox{\boldmath $F$}} }
\newcommand{\bfK}{{\mbox{\boldmath $K$}} }
\newcommand{\bfL}{{\mbox{\boldmath $L$}} }
\newcommand{\bfb}{{\mbox{\boldmath $b$}} }
\newcommand{\bfg}{{\mbox{\boldmath $g$}} }
\newcommand{\bfn}{{\mbox{\boldmath $n$}} }
\newcommand{\bfs}{{\mbox{\boldmath $s$}} }
\newcommand{\cf}{{\it cf.} }
\newcommand{\io}{\int_\Omega}
\newcommand{\1}{\item[({\it i})]}
\newcommand{\2}{\item[({\it ii})]}
\newcommand{\3}{\item[({\it iii})]}
\newcommand{\4}{\item[({\it iv})]}
\newcommand{\5}{\item[({\it v})]}
\newcommand{\6}{\item[({\it vi})]}
\newcommand{\7}{\item[({\it vii})]}
\newcommand{\8}{\item[({\it viii})]}
\newcommand{\9}{\item[({\it xi})]}
\newcommand{\ido}{\int_{\partial\Omega}}
\newcommand{\half}{\mbox{$\frac{1}{2}$}}
\def\parallel{\|}
\def\mid{|}
\def\Bbb R{\real}
\def\hat{\widehat}
\def\tilde{\widetilde}
\def\bar{\overline}
\newcommand{\threehalves}{3\over 2}
\newcommand{\bfPi}{\mbox{\boldmath $\Pi$}}
\newcommand{\bfXi}{\mbox{\boldmath $\Xi$}}
\newcommand{\bfalpha}{\mbox{\boldmath $\alpha$}}
\newcommand{\bfbeta}{\mbox{\boldmath $\beta$}}
\newcommand{\bfgamma}{\mbox{\boldmath $\gamma$}}
\newcommand{\bfdelta}{\mbox{\boldmath $\delta$}}
\newcommand{\bfzeta}{\mbox{\boldmath $\zeta$}}
\newcommand{\bfUpsilon}{\mbox{\boldmath $\Upsilon$}}
\newcommand{\bfGamma}{\mbox{\boldmath $\Gamma$}}
\newcommand{\bfcala}{\mbox{\boldmath ${\cal A}$}}
\newcommand{\bfcalm}{\mbox{\boldmath ${\cal M}$}}
\newcommand{\bfcaln}{\mbox{\boldmath ${\cal N}$}}
\newcommand{\bfcalq}{\mbox{\boldmath ${\cal Q}$}}
\newcommand{\bfcalb}{\mbox{\boldmath ${\cal B}$}}
\newcommand{\bfcalc}{\mbox{\boldmath ${\cal C}$}}
\newcommand{\bfcali}{\mbox{\boldmath ${\cal I}$}}
\newcommand{\bfcalg}{\mbox{\boldmath ${\cal G}$}}
\newcommand{\bfcalh}{\mbox{\boldmath ${\cal H}$}}
\newcommand{\bfcalk}{\mbox{\boldmath ${\cal K}$}}
\newcommand{\bfcalt}{\mbox{\boldmath ${\cal T}$}}
\newcommand{\bfcalx}{\mbox{\boldmath ${\cal X}$}}
\newcommand{\bfcall}{\mbox{\boldmath ${\cal L}$}}
\newcommand{\bfcalf}{\mbox{\boldmath ${\cal F}$}}
\newcommand{\bfcalr}{\mbox{\boldmath ${\cal R}$}}
\newcommand{\bfcals}{\mbox{\boldmath ${\cal S}$}}
\newcommand{\bfcalw}{\mbox{\boldmath ${\cal W}$}}
\newcommand{\bfcalu}{\mbox{\boldmath ${\cal U}$}}
\newcommand{\bfcalv}{\mbox{\boldmath ${\cal V}$}}
\newcommand{\bfcalz}{\mbox{\boldmath ${\cal Z}$}}
\pagenumbering{roman}
\newcommand{\art}[6]{{\I[{\sc #1,}] {#2}, {\it #3}, {\bf #4}, {#5} {[#6]}}}
\newcommand{\ED}{\end{description}}
\newcommand{\I}{\item }
\newcommand{\ra}{\rm a}
\newcommand{\rb}{\rm b}
\newcommand{\rc}{\rm c}
\newcommand{\Hsp}{{\rm I}\!\!\,{\rm R}^n_+}
\newcommand{\Hsn}{{\rm I}\!\!\,{\rm R}^n_-}
\newcommand{\po}[1]{\mbox{$\displaystyle \frac{\mbox{$\partial #1$}}
{\mbox{$\partial x_{1}$}}$}}
\newcommand{\PO}[1]{\mbox{$\displaystyle \frac{\mbox{$\partial #1$}}
{\mbox{$\partial y_{1}$}}$}}
\newcommand{\OP}{\left(\Delta+2\lambda\PO{}\right)}
\newcommand{\op}{\left(\Delta+2\lambda\po{}\right)}
\newcommand{\ft}[1]{
\Frac{1}{(2\pi)^{n/2}}\Int{{\Bbb R}^{n}}{}e^{i{\bf x}\cdot \bfxi}
#1(\xi)d\xi}
\newcommand{\Ft}[1]{
\Frac{1}{2\pi}\Int{{\Bbb R}^{2}}{}e^{i{x}\cdot \xi}
#1(\xi)d\xi}
\newcommand{\Z}{\item[({\it a})]}
\newcommand{\B}{\item[({\it b})]}
\newcommand{\C}{\item[({\it c})]}
\newcommand{\D}{\item[({\it d})]}
\newcommand{\E}{\item[({\it e})]}
\newcommand{\G}{\item[({\it g})]}
\newcommand{\Š}{\`e}
\newcommand{\…}{\`a}
\newcommand{\•}{\`o}
\newcommand{\—}{\`u}
\newcommand{\}{\`{\i}}
\def\tag{\renewcommand{\theequation}}
\newcommand{\Footnote}{~\footnote}
\newcommand{\ie}{{\it i.e.}}
\newcommand{\dist}{\mbox{\rm dist\,}}
\newcommand{\const}{\mbox{\rm const}}
\newcommand{\trace}{\mbox{\rm trace}}
\newcommand{\Bo}{\par\hfill{$\Box$}\par\noindent}
\newcommand{\Nor}[1]{\langle{#1}\rangle_q}
\newcommand{\vs}{\vspace*{.5cm}\par\noindent}
\newcommand{\Vs}{\vspace*{.6cm}\par\noindent}
\newcommand{\Vvs}{\vspace*{.7cm}\par\noindent}
\newcommand{\VVs}{\vspace*{.8cm}\par\noindent}
\newtheorem{definition}{Definition}[section]
\newcommand{\Bd}{\begin{definition}\begin{rm}}
\newcommand{\Ed}{\end{rm}\end{definition}}
\newtheorem{remark}{Remark}[section]
\newcommand{\Br}{\begin{remark}\begin{rm}}
\newcommand{\Er}{\end{rm}\end{remark}}
\newtheorem{proposition}{Proposition}[section]
\newcommand{\Bp}{\begin{proposition}\begin{sl}}
\newcommand{\EP}[1]{\end{sl}\label{proposition:#1}\end{proposition}}
\newcommand{\propref}[1]{{\rm Proposition \ref{proposition:#1}}}
\newcommand{\Bt}{\begin{theorem}\begin{sl}}
\newcommand{\Et}{\end{sl}\end{theorem}}
\newcommand{\Bl}{\begin{lemma}\begin{sl}}
\newcommand{\El}{\end{sl}\end{lemma}}
\newtheorem{theorem}{Theorem}[section]
\newtheorem{lemma}{Lemma}[section]
\newtheorem{corollary}{Corollary}[section]
\newcommand{\eqref}[1]{{\rm (\ref{eq:#1})}}
\newcommand{\Bc}{\begin{corollary}\begin{sl}}
\newcommand{\Ec}{\end{sl}\end{corollary}}
\newcommand{\ET}[1]{\end{sl}\label{theorem:#1}\end{theorem}}
\newcommand{\EDD}[1]{\end{rm}\label{definition:#1}\end{definition}}
\newcommand{\EL}[1]{\end{sl}\label{lemma:#1}\end{lemma}}
\newcommand{\theoref}[1]{{\rm Theorem \ref{theorem:#1}}}
\newcommand{\ER}[1]{\end{rm}\label{remark:#1}\end{remark}}
\newcommand{\EC}[1]{\end{sl}\label{corollary:#1}\end{corollary}}
\newcommand{\remref}[1]{{\rm Remark \ref{remark:#1}}}
\newcommand{\cororef}[1]{{\rm Corollary \ref{corollary:#1}}}
\newcommand{\lemmref}[1]{{\rm Lemma \ref{lemma:#1}}}
\newcommand{\essup}[1]{{\rm ess}\,{{\displaystyle \sup_{\hspace*{-5mm}{#1}}}}}

\renewcommand{\i}{{\rm i}}
\newcommand{\bfsfu}{{\textbf{\textsf{u}}}}
\newcommand{\bfsff}{{\textbf{\textsf{f}}}}
\newcommand{\bfsfF}{{\textbf{\textsf{F}}}}


\pagenumbering{arabic}
\newcommand{\QED}{{\par\hfill$\square$\par}}
\renewcommand{\thefootnote}{(\arabic{footnote})}
\title{Viscous Flow past a Body Translating by Time-Periodic Motion with Zero Average
}
 
\author{ Giovanni P. Galdi 
\thanks{Department of Mechanical Engineering and Materials Science, University of Pittsburgh, PA 15261. 
}}
\date{}
\maketitle
\begin{abstract} We study existence, uniqueness, regularity and asymptotic spatial behavior of  a Navier-Stokes flow past a body moving by a time-periodic translational motion of period $T$, and with zero average. For example, $\mathscr B$ moves in an  oscillating fashion. The flow  is also time-periodic with  same period $T$. However, sufficiently ``far" from the body, the oscillatory component decays faster than the averaged component, so that the flow shows there a distinctive steady-state character. This provides a rigorous  proof of the ``steady streaming" phenomenon.  
\end{abstract}
\renewcommand{\theequation}{\arabic{section}.\arabic{equation}}
\setcounter{section}{0}
\section{Introduction} 
Consider a body, $\mathscr B$, fully immersed in an unbounded Navier-Stokes liquid otherwise at rest,  moving by translational motion with velocity $\bfxi=\bfxi(t)$.  Suppose  $\bfxi$ is time-periodic with period $T$, and that its average over a period of time, $\bar{\bfxi}$, is zero. 
For example,  the direction of $\bfxi$ may be constant, in which case $\mathscr B$ oscillates between two fixed configurations. More generally, the center of mass of $\mathscr B$ moves periodically along a given closed curve, without $\mathscr B$ being able to spin. \par
The question  we would like to address is whether the liquid will execute a corresponding unique  time-periodic regular motion, and what will the flow characteristic be at ``large" spatial distance away from $\mathscr B$.\par 
From the mathematical viewpoint, this question leads us to investigate the same properties for solutions $(\bfu,p)$ to the following set of equations
\be\ba{cc}\smallskip\left.\ba{ll}\medskip
{\partial}_t\bfu-\bfxi(t)\cdot\nabla\bfu+\bfu\cdot\nabla\bfu=\Delta\bfu-\nabla {p}+\bfb\\
\Div\bfu=0\ea\right\}\ \ \mbox{in $\Omega\times (-\infty,\infty)$}\\
\bfu(x,t)=\bfxi(t)\,,\ \ (x,t)\in \partial\Omega\times (-\infty,\infty)\,.
\ea
\eeq{PR}
Here, $\bfu$ and $p$ are velocity and pressure\Footnote{Divided by the constant density of the liquid.} fields of the liquid, respectively, while $\Omega$ is the flow region, namely, the entire space outside $\mathscr B$.\Footnote{For simplicity, we set the coefficient of kinematic viscosity to be 1, since its actual value is entirely irrelevant to our aims.} Moreover, for completeness and also for allowing the special case $\bfxi\equiv\0$, we have included a body force $\bfb=\bfb(x,t)$ which we take to be periodic of the same period $T$.
\par
Despite the very simple formulation, the problem, in its entirety, does not seem to be solvable by the methods currently available, for several reasons that we  explain next. 
\par 
The first contribution to this type of questions --when $\mathscr B$ moves in an unbounded viscous liquid in a {\em non-trivial} time-periodic fashion\footnote{It must be emphasized that if $\mathscr B$ is kept at rest $(\bfxi\equiv\0$ in our case) or is absent ($\bfxi\equiv\0$ and $\Omega=\real^3$), then problems of existence and uniqueness have been successfully addressed and solved, under different assumptions, by a number of authors; see, e.g., \cite{Ma,KoNo,Y,GaSo,KMT} and the review paper \cite{GaKy}.}-- can be found in \cite{GS1}, in the general context where $\mathscr B$ is also allowed to rotate. The tool used there is the so called ``invading domains" technique, based on the Galerkin method coupled with suitable energy estimates. However, by its own nature, such a method is not capable of furnishing enough information on the spatial behavior of the solution ``far" from $\mathscr B$. As a consequence, while the existence of  weak solutions (for ``arbitrary" data) and  strong solutions (for data of restricted ``size")  can be firmly secured, the question of their uniqueness, which actually requires a certain amount of asymptotic spatial ``regularity," is left open and still remains such. For the same reason, the spatial behavior of these solutions at large distances is still not known.
\par 
More recently, two distinct and equally powerful approaches to the study of  time-periodic flow past a body have been independently developed by several authors.\par 
The first one \cite{GaA,Kyed,GaKy,GaKy1}, consists in splitting the velocity field into its averaged component (over a period), $\bar{\bfu}$, and  oscillatory one, ${\bfw}$,  with zero average. The crucial property showed in those articles   is the validity of maximal $L^q$-regularity  for the relevant linearized (time-dependent) problem obeyed by $\bfw$.  As a result, the authors prove that existence and uniqueness of solutions to the full nonlinear problem (for ``small" data), is reduced to show that the  steady-state problem satisfied by $\bar{\bfu}$  is  well-posed  in appropriate homogeneous Sobolev spaces; see \cite{GaA,GaKy1}. Now, for the problem treated here, this theory would work fine if $\bar{\bfxi}\neq \0$ (as showed in \cite{GaA,GaKy1}), thanks to the fact that, in such a case, the linearized steady-state operator is of the Oseen type, for which  well-posedness is a classical result \cite{GaB}. However, our current assumption  requires $\bar{\bfxi}=\0$, and then the pertinent  linearized operator becomes of the Stokes type, for which well-posedness does {\em not} hold \cite{Gnwp}.  
\par 
Another, and entirely different line of attack, traces back to the remarkable paper  \cite{Y}. It is based on a clever duality argument applied to the mild (very weak) formulation of the problem, coupled with {\em sharp} time-decay properties ($L^p-L^q$-estimates) of the evolution operator associated to the relevant linear problem, and of its first spatial derivatives. It must be emphasized that these estimates play a pivotal role for the success of the method. Such an approach,  further refined, generalized and improved by several authors \cite{KMT,NTH,GH,HG}, is particularly effective, because it allows one to establish existence and uniqueness of mild time-periodic (and almost-periodic) solutions when $\mathscr B$ is permitted to translate {and also} rotate, on condition that both translation and rotation vectors be {\em time independent}. However,  its  extension to the time dependent framework is not at all obvious and probably questionable, since  {\em sharp} $L^p-L^q$-estimates  in this more general context  are not necessarily available.  \cite[Theorem 2.2 and Remark 2.1]{H}. 
\par
The method that here we propose and use  is based upon a two-fold strategy. Since, eventually, the nonlinear analysis will be carried out by a contraction mapping argument, it is sufficient to develop this strategy for the relevant linear problem $\mathscr L\!\mathscr P$, say; see \eqref{4.7}. Thus, in the first place,  we establish a number of ``energy estimates" that, once combined with  the ``invading domains" technique of \cite{GS1}, allows us to show existence,  uniqueness and corresponding estimates of time-periodic solutions $(\bfu,p)$ to $\mathscr L\!\mathscr P$ in a very regular function class, provided $\bfxi=\bfxi(t)$ and the ``body force" $\bff$ are sufficiently smooth (see \lemmref{4.2}). Successively, assuming that $\bff$ possesses suitable spatially asymptotic decay properties, we prove that similar properties must hold also for $(\bfu,p)$. This result  --fundamental to the proof of all our main findings--  is obtained as follows. By a classical ``cut-off" argument applied to $\mathscr L\!\mathscr P$, we obtain a similar problem,  $\mathscr L\!\mathscr P_0$, formulated in the whole space $\real^3$; see \eqref{Xm1}--\eqref{frt0}. Furthermore, with the change of coordinates $\bfx\to \bfy:=\bfx-\int_0^t\bfxi(s)ds$, we may absorb the convective term $\bfxi\cdot\nabla\bfu$ in the time derivative, thus reducing the original system of equations in $\mathscr L\!\mathscr P_0$ to a classical Stokes system; see \eqref{Xm3}. By using the basic properties of the fundamental solution associated to the latter, we then show that all solutions to the corresponding Cauchy problem with vanishing initial data must,  along with their first and second spatial derivatives,  decay  algebraically fast at large spatial distances, uniformly in time, with corresponding estimates; see \lemmref{2.3}. The decay is, of course, with respect to the $y$-coordinates. However, just thanks to the fact that $\bfxi$ has zero average, one easily shows that $y$- and $x$-coordinates are ``equivalent" at large distances; see \eqref{as}. Moreover, we prove that the solution to the Cauchy problem must tend, as time goes to infinity, to the time-periodic one of problem $\mathscr L\!\mathscr P_0$, which, in turn, for all $x$ away from the boundary,  coincides with the solution $(\bfu,p)$ to the original problem $\mathscr L\!\mathscr P$. This result, combined with the global regularity of $(\bfu,p)$, finally furnishes the desired uniform spatial decay estimates on the whole domain $\Omega$; see \propref{1}.

With such a complete theory for the linear problem, we can then employ the contraction mapping theorem in a ball of a suitable Banach space, $\mathscr X$, to extend the result to the fully nonlinear case. In this way, in \theoref{m}, we show that if the data $\bfxi$ and $\bfb$ are sufficiently regular and ``small in size," then problem \eqref{PR} possesses one and only one time-periodic solution $(\bfu,p)$ of period $T$ with  $\bfu\in \mathscr X$. In addition, the spatial derivatives of $\bfu$ of order $m=0,1,2$ decay like $|x|^{-m-1}$, uniformly in time. Likewise, $p$ and $\nabla p$ decay as $|x|^{-2}$ and $|x|^{-3}$, respectively, also uniformly in time. 

Our approach also allows us to furnish the   far--field structure of the solution. More precisely, in \theoref{m1}, we prove that $\bfu$ can be decomposed as 
\be
\bfu(x,t)=\bfU(x)+\bfsigma(x)+\bfw(x,t)
. 
\eeq{dec}
where $\bfU$ is the velocity field of a specific steady-state problem (see \lemmref{4.3}),  decaying like $|x|^{-1}$,  $\bfsigma$ is also time independent and decays like $|x|^{-1-\alpha}$, for some $\alpha\in(0,1)$, while $\bfw$ is the oscillatory component of $\bfu$, given by subtracting to $\bfu$ its (time) average,  and decays {\em faster}, like $|x|^{-2}$. The field $\bfU$ is determined up to a (possible other) velocity field, $\bfU_1$, such that $\bfU-\bfU_1$ falls like $|x|^{-1-\delta}$, for some $\delta\in(0,1)$. This analysis shows, in particular, the distinctive steady-state behavior of the far field  solution, thus providing a rigorous formulation  of the steady streaming phenomenon \cite[Chapter XV]{Sch}, \cite{Ri}. In the (less relevant)  case $\bfxi\equiv\0$, we show that $\bfU$ is uniquely determined as the velocity field of a specific Landau solution \cite{KS,KMT}. Moreover, in this situation, we also prove that the oscillatory component $\bfw$ decays even faster, like $|x|^{-3}$, thus sharpening analogous results of \cite{KMT}. 
\par  
The outline of the paper is as follows. Section 2 is dedicated to the linear problem obtained from \eqref{PR} by neglecting the nonlinear term. We prove existence,  uniqueness and asymptotic behavior of corresponding time-periodic solutions. Successively, in Section 3, we combine this findings with the contraction mapping theorem and prove analogous properties for the full nonlinear problem \eqref{PR}, provided $\bfxi$ and $\bfb$ are sufficiently regular and of restricted ``size." In the final Section 4, we give a detailed analysis of the behavior of our solutions at large spatial distances from $\mathscr B$ that shows  the peculiar steady-state character of the flow sufficiently ``far" from $\mathscr B$. 
\setcounter{equation}{0}
\section{Unique Solvability of the Linear Problem}
We begin to collect the main notation used throughout.
The ball in $\real^3$ of radius $R>0$  centered at the origin is indicated by $B_R$, while $B^R$ stands for its complement.
$\Omega$ is the complement of the closure of a bounded  domain $\Omega_0\subset\mathbb R^3$. We shall assume  $\Omega$ of class $C^4$,\footnote{Some of the peripheral results we shall prove require less regularity, but this is irrelevant for our final objective.} and  take the origin of the coordinate system in $\Omega_0$.  We indicate by $B_{R_*}$ a ball containing the closure of $\Omega_0$. For $R\ge R_*$, we set 
$
\Omega_R=\Omega\cap B_R\,,\ \ \Omega^R=\Omega\cap{B^R}$.
Next, for a domain $A\subseteq \real^3$, by $L^q (A)$,  $1\leq q \leq \infty,$  
$W^{m,q}({A}),$ $W_0^{m,q}(A)$, $m \geq 0,$  $(W^{0,q}\equiv W^{0,q}_0\equiv L^q$), we denote usual Lebesgue and Sobolev classes, with corresponding norms $\|.\|_{q,A}$ and $\|.\|_{m,q,A}$.\Footnote{We shall use the same font style to denote scalar, vector and tensor
function spaces.} By the letter $P$ we indicate the (Helmholtz) projector from $L^2(A)$ onto its subspace constituted by solenoidal (vector) function with vanishing normal component, in distributional sense, at $\partial A$.   
 We also set $\int_{A}u\cdot v=\langle u,v\rangle_{A}$.  
$D^{m,2}(A)$ stands for the space of (equivalence classes of) functions $u$ such that
$ 
\sum_{|k|=m}\|D^k u\|_{2,A}<\infty\,.
$ 
Obviously, the latter defines a seminorm in $D^{m,2}(A)$. Also,  by $D_0^{1,2}(A)$ we denote the completion of $C_0^\infty(A)$ in the norm $\|\nabla(\cdot)\|_2$. In the above notation,  the subscript ``$A$" will be omitted, unless confusion arises. A function $u:A\times \real\mapsto \real^3$ is 
{\em $T$-periodic}, $T>0$, if $u(\cdot,t+T)=u(\cdot\,t)$, for a.a. $t\in \real$,
 and we shall denote by $\bar u$ its average over $[0,T]$, namely,
$$
{\bar u}:=\Frac{1}{T}\int_{0}^{T}u(t)dt\,.
$$
Let $B$ be a function space endowed with seminorm $\|\cdot\|_B$. For $r=[1,\infty]$, $T>0$, $L^r(B)$ is the class of functions
$u:(0,T)\rightarrow B$ such that 
$$
\|u\|_{L^r(B)}\equiv\left\{\ba{ll}\smallskip\big( \Int{0}{T}\|u(t)\|_B^r \big)^{\frac 1r}<\infty, \ \ \mbox{if 
$q\in [1,\infty)\,;$}\\   
\essup{t\in[0,T]}\,\|u(t)\|_B <\infty, \ \ \mbox{if $r=\infty.$}
\ea\right.
$$
Likewise, we put
$$
W^{m,r}(B)=\Big\{u\in L^{r}(B): \sum_{k=0}^m\partial_t^ku\in L^{r}(B)\Big\}\,.
$$
Finally, for $A:=\Omega,\real^3$ and $m\ge 1$, we set
$$
[\!]f[\!]_m:=\sup_{x\in A}\,|(1+|x|^m)f(x)|\,,\ \ [\!]f[\!]_{\infty,m}:=\sup_{(x,t)\in A\times (0,T)}\,|(1+|x|^m)f(x,t)|\,.
$$
\par
We now turn to the  main objective of this section that consists in showing existence and uniqueness of $T$-periodic solutions, in appropriate function classes, to the following set of linear equations:
\be\ba{cc}\smallskip\left.\ba{ll}\medskip
{\partial}_t\bfu-\bfxi(t)\cdot\nabla\bfu=\Delta\bfu-\nabla {p}+\bff\\
\Div\bfu=0\ea\right\}\ \ \mbox{in $\Omega\times (0,T)$}\\
\bfu(x,t)=\bfxi(t)\,,\ \ (x,t)\in \partial\Omega\times [0,T]\,,
\ea
\eeq{4.7}
where $\bff=\bff(x,t)$ and $\bfxi=\bfxi(t)$  are suitably prescribed $T$-periodic functions.\smallskip\par To reach this goal, we need a few  preparatory lemmas. 
\Bl
Let $\bfxi \in W^{2,2}(0,T)$ be $T$-periodic. There exists a solenoidal, $T$-periodic function $\tilde{\bfu} \in W^{2,2}(W^{m,q}),$ $m\in \mathbb{N},$ $q\in[1,\infty],$ such that 
$$\ba{ll}\medskip
\tilde{\bfu}(x,t)=\bfxi(t)\,,\ (t,\bfx)\in[0,T]\times\partial\Omega\\ \medskip
\tilde{\bfu}(x,t)=0\,,\ \mbox{for all $t\in [0,T]$, all $|\bfx|\ge\rho$, and some $\rho>R_*$}\,,\\ 
\| \tilde{\bfu} \|_{W^{3,2}(W^{m,q})}\leq C\,\|\bfxi\|_{W^{3,2}(0,T)}
\,,\ea
$$
where $C=C(\Omega,m,q)$. 
\EL{ext}
{\em Proof.} See \cite[Lemma 2.2]{GS1}.\QED
\Bl Let 
$$\bff=\Div \bfcalf\in W^{2,2}(L^{2})\cap L^\infty(D^{1,2})\cap W^{1,2}(D^{2,2}), \ \mbox{with $\bfcalf\in L^{2}(L^{2})$\,,}
$$
and $\bfxi\in W^{3,2}(0,T)$ be prescribed $T$-periodic functions. Then, there exists at least one $T$-periodic 
\be
\bfu\in  W^{2,2}(D^{2,2})\cap W^{1,2}(D^{4,2})\cap L^\infty(L^6)\,,\ \nabla\bfu\in L^{\infty}(W^{2,2})\,,\ \partial_t\bfu\in W^{1,\infty}(W^{1,2})\,,
\eeq{class}
solving \eqref{4.7}
for a corresponding $T$-periodic function 
\be
p\in L^\infty(L^6\cap D^{1,2}\cap D^{2,2})\cap W^{1,2}(D^{3,2})\,
. 
\eeq{class_1}
Moreover, the solution $(\bfu,p)$ satisfies the following estimate
\be\ba{ll}\smallskip
\|\partial_t\bfu\|_{W^{1,\infty}(W^{1,2})}+\|\bfu\|_{L^{\infty}(L^6)}+\|\nabla\bfu\|_{L^\infty(W^{2,2})} +\|p\|_{L^\infty(L^6)}+\|\nabla p\|_{L^\infty(W^{1,2})}\\ \smallskip
\hspace*{2cm}+\|\bfu\|_{W^{2,2}(D^{2,2})}+\|\bfu\|_{W^{1,2}(D^{4,2})}+\Sum{|k|=3}{}\|D^kp\|_{W^{1,2}(L^2)}\\
\hspace*{1cm}\le C\,\big(\|\bff\|_{W^{2,2}(L^2)}+\|\bff\|_{L^{\infty}(D^{1,2})}+\|\bff\|_{W^{1,2}(D^{2,2})}+\|\bfcalf\|_{L^{2}(L^{2})}+\|\bfxi\|_{W^{3,2}(0,T)} \big)
\ea
\eeq{est}
where $C=C(\Omega,T,\xi_0)$, 
for any fixed $\xi_0$ such that $\|\bfxi\|_{W^{2,2}(0,T)}\le \xi_0$. Finally, if $\int_0^T\bfxi(t)dt=0$, the solution is also unique in the class \eqref{class}, \eqref{class_1}.\footnote{See Footnote \ref{foot:5}.}
\EL{4.2}
{\em Proof.} The proof of  existence  is obtained by an argument similar to that employed in \cite[Sections 3 \& 4]{GS1},  that combines the Galerkin method with the ``invading domains" procedure. Specifically,  we write $\bfu=\bfv+\tilde{\bfu}$, with $\tilde{\bfu}$ given in \lemmref{ext}, and begin to consider problem \eqref{4.7} along an increasing, unbounded sequence of (bounded) domains $\{\Omega_{R_k}\}$ with $\cup_{k\in\nat}\Omega_{R_k}=\Omega$, namely,
\be\ba{cc}\smallskip\left.\ba{ll}\medskip
{\partial}_t\bfv_k-\bfxi(t)\cdot\nabla\bfv_k=\Delta\bfv_k-\nabla \tilde{p}_k+\bff+\bff_c\\
\Div\bfv_k=0\ea\right\}\ \ \mbox{in $\Omega_{R_k}\times (0,T)$}\\
\bfv_k(x,t)=\0\,,\ \ (x,t)\in \partial\Omega_{R_k}\times [0,T]\,,
\ea
\eeq{4.8}
where 
$$
\bff_c:= \Delta\tilde{\bfu}-\partial_t\tilde{\bfu}+\bfxi(t)\cdot\nabla\tilde{\bfu}
$$
If we formally dot-multiply \eqref{4.8}$_1$ by $\bfv_k$ and integrate by parts over $\Omega_{R_k}$ we get
\be
\half\ode{}t\|\bfv_k(t)\|_2^2+\|\nabla\bfv_k(t)\|_2^2=\langle \bff+\bff_c,\bfv_k\rangle\le c_0\left(\|\bfcalf\|_2+\|\bff_c\|_{\frac65}\right)\|\nabla\bfv_k\|_2\,, 
\eeq{4.9}
where we have used the assumption on $\bff$ and the Sobolev inequality 
\be
\|\bfz\|_6\le \gamma_0\,\|\nabla\bfz\|_2,\ \ \ \bfz\in D_0^{1,2}(\real^3)\,, 
\eeq{So}
with $\gamma_0$ numerical constant. 
Employing in \eqref{4.9} Cauchy inequality along with Poincar\`e inequality $\|\bfv_k\|_2\le c_{R_k}\|\nabla\bfv_k\|_2$  we get, in particular,
$$
\ode{}t\|\bfv_k(t)\|_2^2+c_{1R_k}\|\bfv_k(t)\|_2^2\le c_2\,\left(\|\bfcalf\|_2^2+\|\bff_c\|_{\frac65}^2\right)\,.
$$
Combining this inequality with Galerkin method one thus shows the existence of a $T$-periodic (distributional) solution $\bfv_k$ to \eqref{4.8} with $\bfv_k\in L^{\infty}(L^2(\Omega_{R_k}))\cap L^2(D_0^{1,2}(\Omega_{R_k}))$ (see \cite[Lemma 3.1]{GS1}).  Furthermore,
\be
\|\nabla\bfv_k\|_{L^{2}(L^{2})}\le c\,\left(\|\bfcalf\|_{L^{2}(L^{2})}+\|\bff_c\|_{L^{2}(L^{\frac65})}\right)\,,
\eeq{4.10}
where the constant $c$ is independent of $R_k$; see \cite[Section 3]{GS1} for technical details. Notice that, by the mean value theorem, from \eqref{4.10} it follows that there is $t_0\in (0,T)$ such that
\be
\|\nabla\bfv_k(t_0)\|_2^2\le \,c_3\left(\|\bfcalf\|_{L^{2}(L^{2})}^2+\|\bff_c\|^2_{L^{2}(L^{\frac65})}\right)\,.
\eeq{mvt}
In order to obtain more regular solutions, we need to show uniform (in $k$)  estimates for $\bfv_k$ in spaces of higher regularity. For this, we formally dot-multiply \eqref{4.8}$_1$ one time by $P\Delta\bfv_k$, a second time by $\partial_t\bfv_k$ and integrate by parts over $\Omega_{R_k}$. We thus show
\be\ba{ll}\medskip
\half\ode{}t\|\nabla\bfv_k(t)\|_2^2+\|P\Delta\bfv_k(t)\|_2^2=\langle\bfxi\cdot\nabla\bfv_k,P\Delta\bfv_k\rangle+\langle\bff+\bff_c,P\Delta\bfv_k(t)\rangle\\
\half\ode{}t\|\nabla\bfv_k(t)\|_2^2+\|\partial_t\bfv_k(t)\|_2^2=\langle\bfxi\cdot\nabla\bfv_k,\partial_t\bfv_k\rangle+\langle\bff+\bff_c,\partial_t\bfv_k(t)\rangle
 
\ea
\eeq{JB}
which, in turn, by Cauchy-Schwarz inequality entails
$$
\ode{}t\|\nabla\bfv_k(t)\|_2^2+c_4(\big(\|\partial_t\bfv_k(t)\|_2^2+\|P\Delta\bfv_k(t)\|_2^2\big)\le c_5\big(\|\bff\|_2^2+\|\bff_c\|_2^2+\|\nabla\bfv_k(t)\|_2^2\big)\,, 
$$
with $c_5=c_5(\xi_0)$.
We now integrate this differential inequality over $[t_0,t]$,  and use the $T$-periodicity property along with \eqref{mvt} and the inequality
\be
\|D^2\bfz\|_{2,\Omega_R}\le c_{\Omega}\,\left(\|P\Delta\bfz\|_{2,\Omega_R}+\|\nabla\bfz\|_{2,\Omega_R}\right)\,,\ \ \bfz\in D^{1,2}(\Omega_R)\cap D^{2,2}(\Omega_R)\,, 
\eeq{hey}
where $c_{\Omega}$ depends only on the regularity of $\Omega$ \cite[Lemma 1]{Hey} but {\em not} on $R$.
One can thus prove that  $\bfv_k\in W^{1,2}(L^2(\Omega_{R_k}))\cap L^\infty(D_0^{1,2}(\Omega_{R_k}))\cap L^2(D^{2,2}(\Omega_{R_k}))$ and satisfies the uniform bound \cite[Lemma 4.1]{GS1} 
\be\ba{rl}\medskip
\|\bfv_k\|_{L^{\infty}(L^6)}+\|\nabla\bfv_k\|_{L^{\infty}(L^2)}+&\|\partial_t\bfv_k\|_{L^{2}(L^2)}+\|D^2\bfv_k\|_{L^{2}(L^2)}\\ \medskip &\le c\,\big(\|\bff\|_{L^{2}(L^2)}+\|\bfcalf\|_{L^{2}(L^2)}+\|\bff_c\|_{L^{2}(L^\frac65)}\big)\\ &
\le C\,\big(\|\bff\|_{L^{2}(L^2)}+\|\bfcalf\|_{L^{2}(L^{2})}+\|\bfxi\|_{W^{2,2}(0,T)} \big)
\,,\ea
\eeq{4.11}
with $C$ independent of $R_k$ and where, in the last step, we used \lemmref{ext}. Next, we take the time derivative of both sides of \eqref{4.8}$_1$, and dot multiply the resulting equation one time by $\partial_t\bfv_k$, a second time by $P\Delta\partial_t\bfv_k$ and integrate  over $\Omega_{R_k}$. We then obtain
\be
\half\ode{}t\|\partial_t\bfv_k(t)\|_2^2+\|\nabla\partial_t\bfv_k(t)\|_2^2=\langle\bfxi'\cdot\nabla\bfv_k,\partial_t\bfv_k\rangle+\langle\partial_t\bff+\partial_t\bff_c,\partial_t\bfv_k(t)\rangle\,,
\eeq{int}
and
\be\ba{rl}\medskip
\half\ode{}t\|\nabla\partial_t\bfv_k(t)\|_2^2+\|P&\!\!\Delta\partial_t\bfv_k(t)\|_2^2\\&=\langle\bfxi'\cdot\nabla\bfv_k,P\Delta\partial_t\bfv_k(t)\rangle+\langle\partial_t\bff+\partial_t\bff_c,P\Delta\partial_t\bfv_k(t)\rangle\,.
\ea\eeq{chi}
From \eqref{4.11} and the mean value theorem we find that there exists at least one $t_1\in (0,T)$ such that
\be
\|\partial_t\bfv_k(t_1)\|_2\le c\,\big(\|\bff\|_{L^{2}(L^2)}+\|\bfcalf\|_{L^{2}(L^{2})}+\|\bfxi\|_{W^{2,2}(0,T)} \big)\,.
\eeq{spa}
Thus, we integrate \eqref{int} over $[t_1,t]$
and use Cauchy-Schwarz inequality, \eqref{spa},  \eqref{4.11} and the $T$-periodicity of $\bfv_k$, to show
\be
\|\partial_t\bfv_k\|_{L^\infty(L^{2})}+\|\nabla\partial_t\bfv_k\|_{L^{2}(L^2)}\le C\,\big(\|\bff\|_{W^{1,2}(L^2)}+\|\bfcalf\|_{L^{2}(L^{2})}+\|\bfxi\|_{W^{2,2}(0,T)} \big)
\eeq{JB0}
Operating in a similar fashion on \eqref{chi}, and also employing \eqref{JB0} and \eqref{hey}, we get
\be\ba{rl}\medskip
\|\nabla\partial_t\bfv_k\|_{L^\infty(L^{2})}+\|&\!\!D^2\partial_t\bfv_k\|_{L^{2}(L^2)}\\ &\le C\,\big(\|\bff\|_{W^{1,2}(L^2)}+\|\bfcalf\|_{L^{2}(L^{2})}+\|\bfxi\|_{W^{2,2}(0,T)} \big)\,.
\ea\eeq{JB01}
Therefore, combining \eqref{4.11}, \eqref{JB0}, and \eqref{JB01} we infer
\be\ba{ll}\medskip
\|\partial_t\bfv_k\|_{L^\infty(W^{1,2})}+\|\bfv_k\|_{L^\infty(L^6)}+\|\nabla\bfv_k\|_{L^\infty(L^2)}+\|D^2\bfv_k\|_{W^{1,2}(L^2)}\\ \medskip
\hspace*{4cm}\le C\,\big(\|\bff\|_{W^{1,2}(L^2)}+\|\bfcalf\|_{L^{2}(L^{2})}+\|\bfxi\|_{W^{2,2}(0,T)} \big)
\ea
\eeq{JB1}
where $C$ is independent of $k$. By an entirely similar argument, it is now straightforward to show  estimate  \eqref{JB1} with $\bfv_k$, $\bff$ and $\bfxi$ replaced by $\partial_t\bfv_k$, $\partial_t\bff$ and $\bfxi^\prime$. To this end, we first differentiate both sides of \eqref{4.8}$_1$ with respect to time, dot-multiply the resulting equation by $\partial^2_t\bfv_k$ and integrate over $\Omega_{R_k}$ to get
\be\ba{rl}\medskip
\half\ode{}t\|\nabla\partial_t\bfv_k(t)\|_2^2&\!+\|\partial_t^2\bfv_k(t)\|_2^2\\
&=\langle\bfxi'\cdot\nabla\bfv_k+\bfxi\cdot\nabla\partial_t\bfv_k,\partial^2_t\bfv_k\rangle+\langle\partial_t\bff+\partial_t\bff_c,\partial^2_t\bfv_k\rangle\,.\ea
\eeq{cu} 
Successively, by differentiating two times both sides of \eqref{4.8}$_1$ with respect to time and dot-multiplying the resulting equation one time by $\partial^2_t\bfv_k$, a second time by $P\Delta\partial^2_t\bfv_k$, and integrating over $\Omega$, we show
\be\ba{rl}\medskip
\half\ode{}t\|\partial_t^2\bfv_k(t)\|_2^2&\!+\|\nabla\partial_t^2\bfv_k(t)\|_2^2\\ 
&=\langle\bfxi^{\prime\prime}\cdot\nabla\bfv_k+2\bfxi^\prime\cdot\nabla\partial_t\bfv_k,\partial_t^2\bfv_k\rangle+\langle\partial^2_t\bff+\partial^2_t\bff_c,\partial^2_t\bfv_k\rangle\,,
\ea
\eeq{tem}
and
\be
\ba{ll}\medskip
\half\ode{}t\|\nabla\partial_t^2\bfv_k(t)\|_2^2+\|P\Delta\partial_t^2\bfv_k(t)\|_2^2\\=\langle\bfxi^{\prime\prime}\cdot\nabla\bfv_k+2\bfxi^\prime\cdot\nabla\partial_t\bfv_k+\bfxi\cdot\nabla\partial_t^2\bfv_k,P\Delta\partial_t^2\bfv_k\rangle+\langle\partial^2_t\bff+\partial^2_t\bff_c,P\Delta\partial_t^2\bfv_k\rangle\,.
\ea
\eeq{muo}
Thus, using \eqref{cu}--\eqref{muo} and following exactly the same procedure as the one leading to
\eqref{JB1}, one can prove
\be\ba{ll}\medskip
\|\partial^2_t\bfv_k\|_{L^\infty(W^{1,2})}+\|\nabla\partial_t\bfv_k\|_{L^\infty(L^2)}+\|D^2\partial_t\bfv_k\|_{W^{1,2}(L^2)}\\ \medskip
\hspace*{4cm}\le C\,\big(\|\bff\|_{W^{2,2}(L^2)}+\|\bfcalf\|_{L^{2}(L^{2})}+\|\bfxi\|_{W^{3,2}(0,T)} \big)
\ea
\eeq{JB_1}
Finally, setting $\bfF_k:=\Delta\bfv_k+\bff+\bff_c$, from \eqref{4.8}$_1$ we get, formally, that $\tilde{p}_k$ obeys for a.a. $t\in [0,T]$ the following Neumann problem\footnote{Note that $\bfxi(t)\cdot\nabla\bfv_k\cdot\bfn|_{\partial\Omega_{R_k}}=0$.}
\be
\Delta \tilde{p}_k=\Div\bfF_k\ \ \mbox{in $\Omega_{R_k}$}\,,\ \ \partial \tilde{p}_k/\partial\bfn|_{\partial\Omega_{R_k}}=\bfF_k\cdot\bfn\,.
\eeq{4.12}  
Therefore, multiplying both sides of the first equation by $\tilde{p}_k$ and integrating by parts over $\Omega_{R_k}$ we easily establish that the pressure field ${p}_k$ associated to $\bfv_k$ satisfies the estimate \cite[Lemma 4.3]{GS1}
\be
\|\nabla \tilde{p}_k\|_2\le c\,\big(\|D^2\bfv_k\|_2+\|\bff\|_2+\|\bff_c\|_2\big)
\eeq{4.13}
with $c$ independent of $k$. We may now let $R_k\to\infty$ and use the uniform estimate \eqref{JB1} and \lemmref{ext}, 
to show the existence of a pair 
$
(\bfu:=\bfv+\tilde{\bfu},\tilde{p})$, with $\bfu$ $T$-periodic, in the class 
\be
\bfu\in  W^{2,\infty}(D^{1,2})\cap W^{2,2}(D^{2,2})\cap L^{\infty}(L^{6})\,,\ \partial_t\bfu\in W^{1,\infty}(L^{2})\,, \ \   \tilde{p}\in L^2(D^{1,2})\,,
\eeq{JB0}
such that
\be\ba{ll}\medskip
\|\partial_t\bfu\|_{W^{1,\infty}(W^{1,2})}+\|\bfu\|_{L^\infty(L^6)}+\|\nabla\bfu\|_{L^\infty(L^2)}+\|D^2\bfu\|_{W^{2,2}(L^2)}\\ 
\hspace*{4cm}\le C\,\big(\|\bff\|_{W^{2,2}(L^2)}+\|\bfcalf\|_{L^{2}(L^{2})}+\|\bfxi\|_{W^{3,2}(0,T)} \big)\,,
\ea
\eeq{JB2}
and which, in addition, solves the original problem \eqref{4.7}. The proof of this convergence property is entirely analogous to that given in \cite[Lemma 3.4 and Section 4]{GS1}, to which we refer for the missing  details. We shall now prove the $T$-periodicity of the pressure field. To this end, we notice that, for a.a. $t\in [0,T]$, by \cite[Theorem II.6.1]{GaB}, there is a function $p_0=p_0(t)$ such that ${p}:=\tilde{p}-p_0$ satisfies
\be 
\|{p}\|_6\le c_0\,\|\nabla p\|_2\,,
\eeq{4.14}
with $c_0$ depending only on $\Omega$. Proceeding as in the proof of \eqref{4.12},  we recognize that ${p}$ must obey (in the sense of distributions) the problem
$$
\Delta{p}=\Div\bfG\,\ \ \mbox{in $\Omega$}\,,\ \ \partial{p}/\partial\bfn|_{\partial\Omega}=\bfG\cdot\bfn\,,
$$
with $\bfG:=\Delta\bfu+\bfxi\cdot\nabla\bfu-\bfxi'+\bff$.
Since ${p}$ satisfies \eqref{4.14} and $\bfG$ is $T$-periodic, we may exploit a classical uniqueness result and conclude that ${p}$ can be time-wise extended to the entire line to become $T$-periodic as well. In order to complete the existence part of the lemma, we  recall some classical properties of  solutions to the Stokes problem:
\be
\ba{cc}\smallskip\left.\ba{ll}\medskip
\Delta\textbf{\textsf {w}}=\nabla {\sf p}+\textbf{\textsf{F}}\\
\Div\textbf{\textsf {w}}=0\ea\right\}\ \ \mbox{in $\Omega$}\\
\textbf{\textsf {w}}(x)=\textbf{\textsf {w}}_\star\,,\ \ x\in \partial\Omega\,.
\ea
\eeq{hm}
In particular,  we get that any distributional solution to \eqref{hm} satisfies the following estimate for $m=0,1,2$, \cite[Lemma V.4.3]{GaB} 
\be\ba{rl}
\Sum{|k|=0}m\big(\|D^{k+2}\textbf{\textsf {w}}\|_2+\|&\!\!D^{k+1} {\sf p}\|_2\big)\\ \medskip&\le C\,\big(\|\textbf{\textsf{F}}\|_{m,2}+\|\textbf{\textsf {w}}_\star\|_{m-1/2,2,\partial\Omega}+\|\textbf{\textsf {w}}\|_{2,\Omega_R}+\|{\sf p}\|_{2,\Omega_R}\big)
\ea
\eeq{hm1}
Let $h\in L^{2}(\Omega_R)$ with $\int_{\Omega_R}h=0$, and let $\bfphi\in W^{1,2}_0(\Omega_R)$ be a solution to the problem $\Div\bfphi=h$ in $\Omega_R$, satisfying $\|\bfphi\|_{1,2}\le c_R\|h\|_{2}$. The existence of such a $\bfphi$ is well known \cite[Theorem III.3.1]{GaB}. Dot-multiplying both sides of \eqref{hm}$_1$ by $\bfphi$ and integrating by parts over $\Omega_R$, we get
$$
\langle \textbf{\textsf{F}},\bfphi\rangle+\langle\nabla\textbf{\textsf {w}},\nabla\bfphi\rangle=\langle{\sf p},\Div\bfphi\rangle=\langle {\sf p},h\rangle\,.
$$ 
From this relation,  the properties of $\bfphi$ and the arbitrariness of $h$, we deduce that ${\sf p}$, modified by a possible addition of a ($T$-periodic) function of time, must obeys the following inequality 
$$
\|{\sf p}\|_{2,\Omega_R}\le c_R\,\big(\|\textbf{\textsf{F}}\|_{2,\Omega_R}+\|\nabla\textbf{\textsf{w}}\|_{2,\Omega_R}\big)\le C_{R}\,(\|\textbf{\textsf{F}}\|_{2,\Omega_R}\|+\|\textbf{\textsf{w}}\|_{2,\Omega_R}\big)+\half\|D^2\textbf{\textsf{w}}\|_2\,,
$$
where, in the last step, we have used Ehrling inequality. As a result, \eqref{hm1}
furnishes
\be\ba{rl}
\Sum{|k|=0}m\big(\|D^{k+2}\textbf{\textsf {w}}\|_2+\|&\!\!D^{k+1} {\sf p}\|_2\big)\\ \medskip&\le C\,\big(\|\textbf{\textsf{F}}\|_{m,2}+\|\textbf{\textsf {w}}_\star\|_{m-1/2,2,\partial\Omega}+\|\textbf{\textsf {w}}\|_{2,\Omega_R}\big)
\ea
\eeq{Li}
We next observe that, for each $t\in [0,T]$,  \eqref{4.7} can be put in the form \eqref{hm} with
$$
\textbf{\textsf{w}}\equiv \bfu\,,\ \  {\sf p}\equiv p\,,\ \ \textbf{\textsf{F}}\equiv \partial_t\bfu+\bfxi\cdot\nabla\bfu-\bff\,,\ \ \textbf{\textsf{w}}_\star\equiv \bfxi\,, 
$$
so that \eqref{Li} leads to
\be
\ba{ll}\medskip
\Sum{|k|=0}m\left(\|D^{k+2}\bfu(t)\|_2+\|D^{k+1} p(t)\|_2\right)\\
\hspace*{1.0cm}\le C_2\,\big(\|\bff(t)\|_{m,2}+|\bfxi(t)|+\|\partial_t\bfu(t)\|_{m,2}+\|\nabla\bfu(t)\|_{m,2}+\|\bfu(t)\|_{2,\Omega_R}\big)\,,\ea
\eeq{Gy}
with $C_2=C_2(\Omega,m,R,\xi_0)$. If we take $m=0$ in \eqref{Gy} and use \eqref{JB2} we then show
\be
\|D^{2}\bfu\|_{L^\infty(L^2)}+\|\nabla p\|_{L^\infty(L^2)}\le C\,\left(\|\bff\|_{W^{2,2}(L^2)}+\|\bfcalf\|_{L^{2}(L^{2})}+\|\bfxi\|_{W^{3,2}(0,T)} \right)\,.
\eeq{m1}
We next take $m=1$ in \eqref{Gy} and employ  \eqref{JB2} and \eqref{m1} to deduce
\be\ba{rl}\medskip
\Sum{|k|=3}{}\|D^{k}\bfu\|_{L^\infty(L^2)}&\!\!+\|D^2 p\|_{L^\infty(L^2)}\\&\le C\,\left(\|\bff\|_{W^{2,2}(L^2)}+\|\bff\|_{L^{\infty}(D^{1,2})}+\|\bfcalf\|_{L^{2}(L^{2})}+\|\bfxi\|_{W^{3,2}(0,T)} \right)\,.\ea
\eeq{m2}
Finally, \eqref{Gy} with $m=2$ in conjunction with  \eqref{JB2} and \eqref{m2} furnishes
\be\ba{rl}\medskip
\Sum{|k|=4}{}\|D^{k}\bfu\|_{L^2(L^2)}+\Sum{|k|=3}{}\|D^k &\!p\|_{L^2(L^2)}\le C\,\big(\|\bff\|_{W^{2,2}(L^2)}+\|\bff\|_{L^{\infty}(D^{1,2})}\\&+\|\bff\|_{L^2(D^{2,2})}+\|\bfcalf\|_{L^{2}(L^{2})}+\|\bfxi\|_{W^{3,2}(0,T)} \big)\,.\ea
\eeq{m3}
We next consider \eqref{hm} with
$$
\textbf{\textsf{w}}\equiv \partial_t\bfu\,,\ \  {\sf p}\equiv \partial_tp\,,\ \ \textbf{\textsf{F}}\equiv \partial_t\big(\partial_t\bfu+\bfxi\cdot\nabla\bfu-\bff\big)\,,\ \ \textbf{\textsf{w}}_\star\equiv \bfxi^\prime\,, 
$$
and take $m=2$ into \eqref{Li}. Again with the help of \eqref{JB2}, we thus deduce  
\be\ba{rl}\medskip
\Sum{|k|=4}{}\|D^{k}\partial_t\bfu\|_{L^2(L^2)}+\Sum{|k|=3}{}\|D^k \partial_tp&\!\|_{L^2(L^2)}\le C\,\big(\|\bff\|_{W^{2,2}(L^2)}+\|\bff\|_{L^{\infty}(D^{1,2})}\\&+\|\bff\|_{W^{1,2}(D^{2,2})}+\|\bfcalf\|_{L^{2}(L^{2})}+\|\bfxi\|_{W^{3,2}(0,T)} \big)\,.\ea
\eeq{Li_1}
In view of \eqref{m1}--\eqref{Li_1}, the proof of the existence property is thus completed. We shall now prove  uniqueness. This amounts to show that $\bfu\equiv\nabla p\equiv\0$ is the only $T$-periodic solution in the class \eqref{class}, \eqref{class_1} to the problem\Footnote{As a matter of fact, going into the details of the proof, it is readily seen that uniqueness of a solution in the class \eqref{class}--\eqref{class_1} holds in a much larger class than that defined by  \eqref{class}--\eqref{class_1}.\label{foot:5}}
\be\ba{cc}\smallskip\left.\ba{ll}\medskip
{\partial}_t\bfu-\bfxi(t)\cdot\nabla\bfu=\Delta\bfu-\nabla {p}\\
\Div\bfu=0\ea\right\}\ \ \mbox{in $\Omega\times (0,T)$}\\
\bfu(x,t)=\0\,,\ \ (x,t)\in \partial\Omega\times [0,T]\,.
\ea
\eeq{47}
To this end, we 
begin to split $\bfu$ as
\be
\bfu=(\bfu-\bar{\bfu})+\bar{\bfu}:=\bfw+\bar{\bfu}\,
. 
\eeq{471}
Since $\bar{\bfw}=0$, by Poincar\'e inequality,  Fubini's theorem and \eqref{class}, we deduce $\bfw\in L^{2}(L^2)$, so that, in particular,
\be
\bfw\in W^{1,2}(L^2)\cap L^2(W^{2,2})\,.
\eeq{2.15}
From classical embedding theorems (e.g. \cite[Theorem 2.1]{Sol}) and \eqref{2.15} we deduce
\be
\bfw\in L^\infty(L^2)\cap L^s(L^6)\,,\ \ \mbox{all $s\in [2,\infty)$}\,.
\eeq{2.16}
We next observe that from \eqref{47} it follows that $p$ obeys the following Neumann problem for a.a. $t\in [0,T]$
\be
\Delta p=0\ \ \mbox{in $\Omega$}\,,\ \ \pde p\bfn=-\curl\curl\bfu\cdot\bfn\ \ \mbox{at $\partial\Omega$},
\eeq{2.17}
where we used the identity $\Delta\bfu=-\curl\curl\bfu.$ We may modify $p$ by adding to it a suitable $T$-periodic function of time, in such a way that the redefined pressure field, that we continue to denote by $p$, satisfies \eqref{4.14}. Thus, on the one hand, by the mean value theorem, \eqref{class_1}, \eqref{4.14} and smoothness properties of harmonic functions we obtain, in particular,
\be
p\in L^2(C^1(\bar{\Omega_R})\,,\ \ \mbox{for all $R\ge R_*$}\,.
\eeq{2.18}
On the other hand, observing that, by Stokes theorem and \eqref{2.17}, 
$$
0=-\int_{\partial\Omega}\curl\curl\bfu\cdot\bfn=\int_{\partial\Omega}\pde p\bfn=\int_{\partial B_R}\pde p\bfn=0\,,\ \ \mbox{for all $R\ge R_*$}\,,$$ 
from \eqref{4.14} and well-known results on Laplace equation on exterior domains (e.g. \cite[Exercise V.3.6]{GaB}) we find for a.a. $t\in [0,T]$
\be
p(x,t)=\int_{\partial B_R}[(\mathcal E(x-y)-\cale(x))\pde p \bfn(y,t)-p(y,t)\pde\cale \bfn(x-y)]d\sigma_y\,,\ \ |x|\ge 2R\,,
\eeq{2.19}
where $\mathcal E=\cale(z)$ is the Laplace fundamental solution. Since 
\be|\nabla\cale(z)|
\le c\,|z|^{-2}\,, \ \ |z|\neq 0\,,
\eeq{GE}
from \eqref{2.18} and \eqref{2.19} it follows that
\be
p\in L^2(L^r(\Omega^{2R}))\,,\ \ \mbox{all $r> 3/2$}\,.
\eeq{2.20}
Let $\psi_R=\psi_R(x)$ be a smooth cut-off function that is 1 for $|x|\le 2R$,  is 0 for $|x|\ge 3R$ and $|\nabla\psi_R|\le C\,R^{-1}$, with $C$ independent of $R$. Clearly,
\be
\nabla\psi_R\in L^3(\Omega)\,.
\eeq{2.21}
We dot-multiply both sides of \eqref{47}$_1$ by $\psi_R\bfu$, and integrate by parts over $\Omega\times(0,T)$. Noticing that $\bfu\in L^2(L^2(\Omega\rho))$, all $\rho\ge R_*$, and using $T$-periodicity we thus show 
\be\ba{rl}\medskip
\Int0T\Int\Omega{}\psi_R\,|\nabla\bfu|^2&\!=-\half\Int0T\Int{\Omega_{2R,3R}}{}\nabla\psi_R\cdot\bfxi(t)|\bfu|^2+\Int0T\Int{\Omega_{2R,3R}}{} p\,\nabla\psi_R\cdot\bfu\\ :&\!= -\half I_{1R}+I_{2R}\,.
\ea\eeq{2.22}
From H\"older inequality and \eqref{class}
$$
|I_{2R}|\le \sup_{t\in [0,T]}\|\bfu(t)\|_6\|\nabla\psi_R\|_3\int_0^T\|p(t)\|_{2,\Omega^{2R}}\,, 
$$
which, by \eqref{2.21}, entails
\be
\lim_{R\to\infty}|I_{2R}|=0\,.
\eeq{2.23}
Furthermore, employing \eqref{471} and Fubini's theorem, we show
$$\ba{rl}\medskip
I_{1R}&\!=\Int{\Omega_{2R,3R}}{}\nabla\psi_R\cdot\Int0T\bfxi(t)\,(|\bar{\bfu}|^2+|\bfw|^2+2\bar{\bfu}\cdot\bfw)\\ \medskip&=\Int0T\Int{\Omega_{2R,3R}}{}\nabla\psi_R\cdot\bfxi(t)\,(|\bfw|^2+2\bar{\bfu}\cdot\bfw)\\&:=I_{1R}^1+I_{1R}^2\,,
\ea$$
where we have used the assumption $\bar{\bfxi}=\0$. Again by H\"older inequality,
 and the properties of $\psi_R$
$$
|I_{1R}^1|\le c\,\|\bfxi\|_{W^{1,2}(0,T)}\,R^{-1} \big(\int_0^T\|\bfw\|_{2,\Omega^{2R}}^2\big)^{\frac12}\,,
$$
which, 
by \eqref{2.16}, implies
\be
\lim_{R\to\infty}|I_{2R}^1|=0\,.
\eeq{2.24}
Finally, by using one more time H\"older inequality, we infer
$$
|I_{1R}^2|\le c\,\|\bfxi\|_{W^{1,2}(0,T)}\|\nabla\psi_R\|_3\,\|\bfw\|_{L^\infty(L^2)}\big(\int_0^T\|\bar{\bfu}\|_{6,\Omega^{2R}}^6
\big)^\frac16\,,
$$
and so from the latter, \eqref{2.16} and \eqref{class} we obtain
\be
\lim_{R\to\infty}|I_{2R}^2|=0\,.
\eeq{2.25}
Uniqueness then follows by letting $R\to\infty$ in \eqref{2.22} and using \eqref{2.23}--\eqref{2.25}. The lemma is completely proved.
\par\hfill$\square$
\Bl
Let $\bfG$ be a second-order tensor field in $\mathbb{R}^3 \times (0,\infty)$ such that
$$
[\!]\bfG[\!]_{\infty,2} +\sum_{|k|=0}^1 [\!]D^k(\nabla \cdot \bfG)[\!]_{\infty,|k|+3} < \infty\,.
$$
Then, the Cauchy problem 
\be\ba{cc}\medskip\left.\ba{ll}\medskip
\partial_t\bfv=\Delta\bfv-\nabla\phi +\Div\bfG\\
\Div\bfv=0\ea\right\}\ \ \mbox{in $\real\times(0,\infty)$}\\
\bfv(x,0)=\0\,,\ \ x\in\real^3
\ea
\eeq{Cau}
has one and only one solution such that for all $\tau>0$,
\begin{equation}
(\bfv,\phi) \in [W^{1,2}(0,\tau;L^2(\real^3))\cap L^2(0,\tau;W^{2,2}(\real^3))]\times L^2(0,\tau;D_0^{1,2}(\real^3)).
\label{estimunn1}
\end{equation}
Moreover, 
$$
\sum_{|k|=0}^2 [\!]D^k\bfv [\!]_{\infty,|k|+1}+
\sum_{|k|=0}^1 [\!]D^k\phi[\!]_{\infty,|k|+2}<\infty\,,
$$
and the following inequality holds:
\begin{equation}
\sum_{|k|=0}^2 [\!]D^k\bfv [\!]_{\infty,|k|+1}+
\sum_{|k|=0}^1 [\!]D^k\phi[\!]_{\infty,|k|+2} \leq C\,\Big([\!]\bfG[\!]_{\infty,2}+\Sum{|k|=0}{1}[\!]D^k\nabla\cdot\bfG[\!]_{\infty,|k|+3}\Big)\,,
\label{estimunn2}
\end{equation}
with $C$ a (positive) numerical constant. 
\EL{2.3}
{\em Proof.} The existence of a unique solution in the class (\ref{estimunn1}) is a classical result (e.g. \cite[Theorem VIII.4.1]{GaB}). Moreover, the velocity field $\bfv$ admits the following integral representation
\be
v_i(x,t)=\int_0^t\int_{\real^3}\Gamma_{ih}(x-y,t)\partial_{j}G_{jh}(y,t-s)dyds\,, \ \ i=1,2,3\,,
\eeq{rep}
where $\bfGamma=\bfGamma(\chi,\rho)$ is the Oseen fundamental solution to the Stokes problem
 (see \cite[Theorem VIII.4.2]{GaB}) for which, in particular, the following estimates hold:
\be
\Gamma_k(\chi):=\int_0^\infty|D^k\bfGamma(\chi,t)|dt\le C\,|\chi|^{-m};\,\ |k|=m\in \{1,2,3\}\,, \, \ \bfchi\neq\0; 
\eeq{OE}
see \cite[Lemma VIII.3.3 and Exercise VIII.3.1]{GaB}. Using \eqref{rep} and \eqref{OE} along with the assumption on $\bfG$, one can then establish the stated pointwise estimate on $\bfv$  \cite[Theorem VIII.4.4]{GaB}. We shall now prove the claimed property for $\nabla\bfv$. To this end, let $R=\half |x|>1$. From \eqref{rep} we thus deduce 
\be\ba{ll}\medskip
\partial_kv_i(x,t)=\displaystyle{\int_0^t\int_{B_R}}\partial_j\partial_k\Gamma_{ih}(x-y,s)G_{jh}(y,t-s)dyds\\ \medskip
\hspace*{3.5cm}+\displaystyle{\int_0^t\int_{\partial B_R}}\partial_k\Gamma_{ih}(x-y,s)G_{jh}(y,t-s)n_jd\sigma_yds\\ \medskip
\hspace*{4cm}+\displaystyle{\int_0^t\int_{B^R}}\partial_k\Gamma_{ih}(x-y,s)\partial_jG_{jh}(y,t-s)dyds\\
\hspace*{4cm}
:=I_1+I_2+I_3\,.
\ea
\eeq{split}
From \eqref{OE}, the assumption, and the fact that $|x-y|\ge R$, $y\in B_R$,  it follows
\be\ba{rl}\medskip
|I_1|\le &\!C_1\,[\!]\bfG[\!]_{\infty,2}\Int{B_R}{}\Frac{dy}{|x-y|^3(1+|y|^2)}\le 2C_1\,[\!]\bfG[\!]_{\infty,2}\,|x|^{-3}\Int{B_R}{}\Frac{dy}{(1+|y|^2)}\\
\le&\! C_2\,[\!]\bfG[\!]_{\infty,2}\, |x|^{-2}\,.
\ea
\eeq{A1}
By the same token, again using \eqref{OE}, we get
\be
|I_2|\le C_1\,[\!]\bfG[\!]_{\infty,2}\Int{\partial B_R}{}\Frac{d\sigma_y}{|x-y|^2(1+|y|^2)}\le C_2\,[\!]\bfG[\!]_{\infty,2}\, |x|^{-2}\,.
\eeq{A2}
Furthermore, 
$$\ba{rl}\medskip
|I_3|&\le C_1\,[\!]\Div\bfG(t)[\!]_{\infty,3}\, \Int{B^R}{}\Frac{dy}{|x-y|^2(1+|y|^3)}\\
&\le \Frac{C_1}{2R}\,[\!]\Div\bfG[\!]_{\infty,3}\, \Int{\real^3}{}\Frac{dy}{|x-y|^2|y|^2}\,.
\ea$$
As a result, from a well-known theorem on convolutions \cite[Lemma II.9.2]{GaB} applied to the last integral, we infer
\be
|I_3|\le C_2\,[\!]\Div\bfG[\!]_{\infty,3}\,|x|^{-2}\,.
\eeq{A3}
Finally, if $|x|\le 2$, from \eqref{rep} and \eqref{OE} we deduce
$$\ba{rl}\medskip
|\nabla\bfv(x,t)|&\le C_1\,[\!]\Div\bfG[\!]_{\infty,3}\big(\Int{|x-y|\le 5}{}\Frac{dy}{|x-y|^2}+\Int{B^3}{}\Frac{dy}{|y|^2(1+|y|^3)}\big)\\
&\le C_2\,[\!]\Div\bfG[\!]_{\infty,3}\,.
\ea$$
The latter, 
combined with \eqref{rep}--\eqref{A3} thus proves the desired property for $\nabla\bfv$. By the same token we get the estimate for $D^2\bfv$. Actually, from \eqref{rep} we show by a double integration by parts
\be
\ba{ll}\medskip
\partial_l\partial_kv_i(x,t)=\displaystyle{\int_0^t\int_{B_R}}\partial_l\partial_j\partial_k\Gamma_{ih}(x-y,s)G_{jh}(y,t-s)dyds\\ \medskip
\hspace*{3.5cm}+\displaystyle{\int_0^t\int_{\partial B_R}}\partial_l\partial_k\Gamma_{ih}(x-y,s)G_{jh}(y,t-s)n_jd\sigma_yds\\ \medskip
\hspace*{3.5cm}+\displaystyle{\int_0^t\int_{\partial B_R}}\partial_k\Gamma_{ih}(x-y,s)\partial_jG_{jh}(y,t-s)n_ld\sigma_yds\\ \medskip
\hspace*{4cm}+\displaystyle{\int_0^t\int_{B^R}}\partial_k\Gamma_{ih}(x-y,s)\partial_l\partial_jG_{jh}(y,t-s)dyds\\
\hspace*{4cm}
:={\sf T}_1+{\sf T}_2+{\sf T}_3+{\sf T}_4\,.
\ea
\eeq{fs}
Thus, employing \eqref{OE} and  \cite[Lemma II.9.2]{GaB}, we easily show
\be
|{\sf T}_1|\le C_1\,[\!]\bfG[\!]_{\infty,2}R^{-2}\Int{\real^3}{}\frac{dy}{|x-y|^2|y|^2}\le C_2\,\,[\!]\bfG[\!]_{\infty,2}|x|^{-3}\,,
\eeq{fs1}
and, likewise,
\be
|{\sf T}_4|\le C_3\,[\!]\nabla(\Div\bfG)[\!]_{\infty,4}R^{-2}\Int{\real^3}{}\frac{dy}{|x-y|^2|y|^2}\le C_4\,\,[\!]\nabla(\Div\bfG)[\!]_{\infty,4}|x|^{-3}\,.
\eeq{fs2}
Moreover,
\be
|{\sf T}_2|\le C_5\,[\!]\bfG[\!]_{\infty,2}R^{-3}\Int{\partial B_R}{}\frac{d\sigma_y}{(1+|y|^2)}\le C_6\,\,[\!]\bfG[\!]_{\infty,2}|x|^{-3}\,,
\eeq{fs3}
and
\be
|{\sf T}_3|\le C_7\,[\!]\Div\bfG[\!]_{\infty,3}R^{-3}\Int{\partial B_R}{}\frac{d\sigma_y}{(1+|y|^2)}\le C_8\,\,[\!]\Div\bfG[\!]_{\infty,3}|x|^{-3}\,.
\eeq{fs4}
If $|x|\le2$, as in the analogous estimate for $\nabla\bfv$, we show
$$\ba{rl}\medskip
|D^2\bfv(x,t)|&\le C_9\,[\!]\nabla(\Div\bfG)[\!]_{\infty,4}\big(\Int{|x-y|\le 5}{}\Frac{dy}{|x-y|^2}+\Int{B^3}{}\Frac{dy}{|y|^2(1+|y|^4)}\big)\\
&\le C_{10}\,[\!]\nabla(\Div\bfG)[\!]_{\infty,4}\,.
\ea$$
As a result, the claimed estimate for $D^2\bfv$ follows from the latter and \eqref{fs}--\eqref{fs4}.
The estimates for $\phi$ and $\nabla \phi$  are obtained in an entirely similar fashion. In fact, this is a consequence of the following representation, valid  for a.a. $t\in [0,T]$,
\be
\phi(x,t)=-
\int_{\real^3}\partial_j\cale(x-y)\partial_iG_{ij}(y,t)dy\,,
\eeq{p0}
and of the fact that $D^k\cale(\chi)$ satisfies exactly the same  properties as $\Gamma_k(\chi)$ in \eqref{OE}. We therefore shall omit the proof of these estimates, leaving it to the reader as an exercise.  The lemma is proved.\nopagebreak \QED
\smallskip\par
We are now in a position to show the main result of this section. Precisely, we have the following.
\Bp Let $\bfcalf$ and $\bfxi$ be prescribed $T$-periodic functions such that
$$\ba{ll}\medskip
\bff:=\Div\bfcalf\in W^{2,2}(L^{2})\cap W^{1,2}(D^{2,2})\,,\ 
[\!]\bfcalf[\!]_{\infty,2}+\Sum{|k|=0}{1} [\!]D^k\bff[\!]_{\infty,|k|+3}<\infty\,;
\\ \bfxi\in W^{3,2}(0,T)\,,\ \Int0T\bfxi(t)dt=\0\,.
\ea$$
Then, problem \eqref{4.7} has one and only one solution $(\bfu,p)$ in the class \eqref{class}, \eqref{class_1}, which satisfies the estimate
\be\ba{ll}\smallskip
\|\partial_t\bfu\|_{W^{1,\infty}(W^{1,2})}+\|\bfu\|_{L^{\infty}(L^6)}+\|\nabla\bfu\|_{L^\infty(W^{2,2})} +\|p\|_{L^\infty(L^6)}+\|\nabla p\|_{L^\infty(W^{1,2})}\\ \smallskip
\hspace*{2cm}+\|\bfu\|_{W^{2,2}(D^{2,2})}+\|\bfu\|_{W^{1,2}(D^{4,2})}+\Sum{|k|=3}{}\|D^kp\|_{W^{1,2}(L^2)}\smallskip\\
\hspace*{.5cm}\le C_1\,\big(\|\bff\|_{W^{2,2}(L^2)}+\|\partial_t\bff\|_{L^2(D^{2,2})}+\Sum{|k|=0}{1} [\!]D^k\bff[\!]_{\infty,|k|+3}+[\!]\bfcalf[\!]_{\infty,2}+\|\bfxi\|_{W^{3,2}(0,T)} \big)\\
\hspace*{.5cm}:=C_1\,\mathscr D_1\,.
\ea
\eeq{estik}
In addition, if, for some $\rho\ge R_\star$ $\nabla\bfcalf\in L^\infty(L^\infty(\Omega_{2\rho}))$, then $\sum_{|k|=0}^2[\!]D^k\bfu[\!]_{\infty,|k|+1}$, $\sum_{|k|=0}^1[\!]p[\!]_{\infty,|k|+2}<\infty$, and we have
\be
\Sum{|k|=0}{2}[\!]D^k\bfu[\!]_{\infty,|k|+1}+\Sum{|k|=0}{1}[\!]p[\!]_{\infty,|k|+2}
\le C_2\, \big(\mathscr D+\|\nabla\bfcalf\|_{L^\infty(L^\infty(\Omega_{2\rho})}\big):= C_2\,\mathscr D_2
\,,
\eeq{estikz}
where $C_i=C_i(\Omega,T,\xi_0)$, $i=1,2$.\footnote{Recall that $\xi_0$ is defined in \lemmref{4.2}.}
\EP{1}
{\em Proof.} We begin to observe that, obviously,
\be
\|\bfcalf\|_{L^2(L^2)}\le C\, [\!]\bfcalf[\!]_{\infty,2}\,,\ \ \|\bff\|_{L^2(D^{2,2})}\le C\,\Sum{|k|=0}{1} [\!]D^k\bff[\!]_{\infty,|k|+3}
\,.
\eeq{ha0}
Therefore, under  the given assumptions, the existence and uniqueness of  a solution $(\bfu,p)$ in the class \eqref{class}, \eqref{class_1} satisfying \eqref{estik} is ensured by \lemmref{4.2}.
In order to complete the proof of the proposition, it remains to show the pointwise properties of $\bfu$ and $p$,  along with the corresponding estimates. To this end, for a fixed $\rho\ge R_*$, let $\psi=\psi(x)$ be a smooth ``cut-off" function such that $\psi(x)=0$ for $|x|\le\rho$, $\psi(x)=1$ for $|x|\ge2\rho$, and set $\bfw:=\psi\bfu$, ${\sf p}:=\psi p$. From \eqref{4.7} we thus infer that $(\bfw,{\sf p})$ obeys the following problem
\be\left.\ba{ll}\medskip
{\partial}_t\bfw-\bfxi(t)\cdot\nabla\bfw=\Delta\bfw-\nabla {\sf p}+\Div\bfg+\bfg_c\\
\Div\bfw=h\ea\right\}\ \ \mbox{in $\real^3\times (0,T)$}
\,,
\eeq{Xm1}
where 
\be\ba{ll}\medskip
\bfg:=\psi\,\bfcalf\,,\ \ h:=\nabla\psi\cdot\bfu\\
\bfg_c:=-\nabla\psi\cdot\bfcalf-\bfxi(t)\cdot\nabla\psi\,\bfu-\Delta\psi\,\bfu-2\nabla\psi\cdot\nabla\bfu+p\,\nabla\psi\,. 
\ea
\eeq{frt0}
We next observe that, by classical embedding theorems, 
$$\ba{ll}\medskip
\|\bfu\|_{L^\infty(L^\infty)}\le C\,\big(\|\bfu\|_{L^\infty(L^6)}+\|D^2\bfu\|_{L^\infty(L^2)}\big)\,, \\ \medskip  \|p\|_{L^\infty(L^\infty)}\le C\,\big(\|p\|_{L^\infty(L^6)}+\|D^2p\|_{L^\infty(L^2)}\big)\\ \medskip
\|\nabla\bfu\|_{L^\infty(L^\infty)}\le C\|\nabla\bfu\|_{L^\infty(W^{2,2})}\,,\\ \medskip 
\|\nabla p\|_{L^\infty(L^\infty)}\le C\big(\|\nabla p\|_{L^\infty(L^{2})}+\Sum{|k|=3}{}\|D^kp\|_{W^{1,2}(L^2)}\big)\,,
\\
\|D^2\bfu\|_{L^\infty(L^\infty)}\le C\,\big(\|D^2\bfu\|_{L^\infty(L^2)}+\Sum{|k|=4}{}\|D^k\bfu\|_{W^{1,2}(L^2)}\big)\,.\ea
$$
Therefore, from the latter and \eqref{estik} we get
\be
\Sum{|k|=0}{2}\|D^k\bfu\|_{L^\infty(L^\infty)}+
\Sum{|k|=0}{1}\|D^kp\|_{L^\infty(L^\infty)}
\le C\,\mathscr D_1
\,.
\eeq{Li_2}
We also notice we have
\be\bfg_c=\Div\bfcalh\eeq{0} 
with
\be
[\!] \bfcalh [\!]_{\infty,2}+\Sum{|k|=0}{1}[\!]D^k\bfg_c[\!]_{\infty,|k|+3}\le c\,\mathscr D_2 
\eeq{frt}
where $c=c(\Omega,T,\xi_0)$. In fact, let
$$
\bfcalh(x,t)=\int_{\real^3}\nabla\cale(x-y)\cdot\bfg_c(y,t)dy\,,
$$
where, we recall, $\mathcal E$ is the Laplace fundamental solution. Clearly, $\Div\bfcalh=\bfg_c$ and, by  \eqref{frt0}$_2$, \eqref{Li_2},  $T$-periodicity, and the fact that the support
of $\bfg_c$ is contained in $B_{2\rho}$, it follows at once
\be
\Sum{|k|=0}{1}[\!]D^k\bfg_c[\!]_{\infty,|k|+3}\le C\,\mathscr D_2\,.
\eeq{mm}
Moreover,  from \eqref{GE} 
we find, for a.a. $t\ge 0$
\be
\sup_{|x|\ge 4\rho}|\bfcalh(x,t)|\,|x|^2\le C_1\,\|\bfg_c(t)\|_1
\eeq{ha1}
with $C_1=C_1(\rho)$. Also, from classical results for convolutions with weakly singular integrals (e.g. \cite[Theorem II.11.2]{GaB}), we have
\be
\sup_{|x|\le 4\rho}|\bfcalh(x,t)|\le C_2\,\|\bfg_c(t)\|_q\,,\ \ q>3\,,
\eeq{ha2}
with $C_2=C_2(\rho)$. Consequently, \eqref{frt} follows again from \eqref{frt0}$_2$, \eqref{Li_2}, \eqref{mm}--\eqref{ha2} and $T$-periodicity. Let
\be
\bfV(x,t)=\int_{\real^3}\nabla\cale (x-y)\,h(x,t)
\eeq{ha3}
and write
$$
\bfw(x,t)=\bfw_1(x,t)+\bfV(x,t)\,.
$$
Notice that $\bfV$ is $T$-periodic and, as a result, so is $\bfw_1$. Moreover, using also Sobolev and  Calderon-Zygmund theorems and that $h$ has  bounded (spatial) support, we easily show that $\bfV$ is in the functional class defined  in \eqref{class}. Thus, from \eqref{Xm1}, and taking into account  \eqref{ha0} we deduce that $\bfw_1$ is a $T$-periodic solution in the class \eqref{class} to the following problem
\be\left.\ba{ll}\medskip
{\partial}_t\bfw_1-\bfxi(t)\cdot\nabla\bfw_1=\Delta\bfw_1-\nabla {\sf P}+\Div\bfcalg\\
\Div\bfw_1=0\ea\right\}\ \ \mbox{in $\Omega\times (0,T)$}
\,,
\eeq{Xm2}
with
\be\ba{ll}\medskip
{\sf P}:={\sf p}+\Int{\real^3}{}\cale(x-y)[\partial_th(y,t)-\bfxi\cdot \nabla h(y,t)-\Delta h(y,t)]dy:={\sf p}+\tilde{\sf p}\,,\\
{\bfcalg}:=\bfg+\bfcalh\,.
\ea
\eeq{sa1}
Observe that, by assumption,  \eqref{frt0}$_1$, \eqref{frt} and $T$-periodicity one has
\be 
[\!]\bfcalg[\!]_{\infty,2}
+\Sum{|k|=0}{1}[\!]D^k\Div\bfcalg[\!]_{\infty,|k|+3}
\le C\,\mathscr D_2\,.
\eeq{bfrt}
We now introduce the following change of coordinates
\be
\bfy=\bfx-\bfx_0(t)\,,\ \ \bfx_0(t):=\int_0^t\bfxi(s)ds\,.
\eeq{coo}
Since $\int_0^T\bfxi(t)dt=\0$, this along with the $T$-periodicity of $\bfxi$ implies that $\bfx_0(t)$ is $T$-periodic as well, and also the existence of a constant $M>0$ such that
\be
|\bfx_0(t)|\le M\,.
\eeq{M}
In fact, by integrating over $[0,t]$ both sides of the Fourier series for $\bfxi$:
$$
\bfxi(t)=\Sum{|k|\ge1}{}\bfxi_k\,{\rm e}^{\frac{2\pi}T\,{\rm i}\,k\,t}\,,
$$
we infer at once that $\bfx_0$ is $T$-periodic. Moreover,
$$\ba{rl}\medskip
\big|\Int0t\bfxi(s)ds\big|\le \Frac T{\pi}\Sum{|k|\ge 1}{}\frac1{|k|}|\bfxi_k|&\!\le \Frac T{\pi}\big(\Sum{|k|\ge 1}{}\frac1{|k|^2}\big)^{\frac12}\big(\Sum{|k|\ge 1}{}|\bfxi_k|^2\big)^{\frac12}
\\
&\le C\,\Frac{T^{\frac12}}{\pi}\big(\Int0T|\bfxi(t)|^2\big)^{\frac12}:= M\,.
\ea$$
Notice that from \eqref{coo} and \eqref{M} it follows that
\be 
|\bfy|-M\le|\bfx|\le |\bfy|+M\,. 
\eeq{as}
Setting
\be\ba{ll}\medskip
\bfW_1(y,t):=\bfw_1(y+x_0(t),t)\,,\ \ \ \Pi(y,t):={\sf P}(y+x_0(t),t)\,,\\  \bfG(y,t):=\bfcalg(y+x_0(t),t)\,,
\ea\eeq{def}
from \eqref{coo} and \eqref{Xm2} it follows that $(\bfW_1,\Pi)$ satisfies the following Cauchy problem
\be\ba{cc}\medskip\left.\ba{ll}\medskip
{\partial}_t\bfW_1=\Delta\bfW_1-\nabla {\Pi}+\Div\bfG\\
\Div\bfW_1=0\ea\right\}\ \ \mbox{in $\real^3\times (0,T)$}
\,,\\
\bfW_1(y,0)=\bfw_1(x,0)\equiv \psi(x)\,\bfu(x,0)-\bfV(x,0)\,.
\ea
\eeq{Xm3}
In view of \eqref{coo}, \eqref{def} and  \eqref{as}, we have
\be
|\bfG(y,t)|\,(|\bfy|^2+1)=|\bfcalg(x,t)|\,(|\bfy|^2+1)\le |\bfcalg(x,t)|\,[(|\bfx|+M)^2+1)\le c_1\,[\!]\bfcalg[\!]_{\infty,2}\,,
\eeq{A}
and, likewise,
\be
\Sum{|k|=0}{1}[\!]D^k\Div\bfG[\!]_{\infty,|k|+3}\le c_2\,\Sum{|k|=0}{1}[\!]D^k\Div\bfcalg[\!]_{\infty,|k|+3}\,.
\eeq{B}
As a result, by \eqref{bfrt}, the tensor field $\bfG$ satisfies the assumptions of \lemmref{2.3}. Set 
\be\bfU:=\bfW_1-\bfv\,,\ \ {\sf Q}:=\Pi-\phi\,,
\eeq{saz} 
with $(\bfv,\phi)$ solution given in that lemma. From \eqref{Cau} and \eqref{Xm3}, we then have that $(\bfU,{\sf Q})$ satisfies: 
\be\ba{cc}\medskip\left.\ba{ll}\medskip
{\partial}_t\bfU=\Delta\bfU-\nabla {\sf Q}\\
\Div\bfU=0\ea\right\}\ \ \mbox{in $\real^3\times (0,\infty)$}
\,,\\
\bfU(y,0)= \psi(x)\,\bfu(x,0)-\bfV(x,0)\,.
\ea
\eeq{Xm4}
Since both $\bfu$ and $\bfV$ are in the function class defined by \eqref{class}, we have, in particular,
$$
\bfU(y,0)\in L^6(\real^3)\,,
$$
so that, by classical results on the Cauchy problem for Stokes equations (e.g., \cite[Theorem VIII.4.3]{GaB}) we infer
$$
\lim_{t\to\infty}\left(\Sum{|\ell|=1}{3}\|D^\ell\bfU(t)\|_6+\Sum{|\ell|=0}{2}\|D^\ell{\sf Q}(t)\|_6+\|\bfU(t)\|_r\right)=0\,,\ \ r>6\,.
$$
which, in turn, by embedding, implies
\be
\lim_{t\to\infty}\Big(\Sum{|\ell|=0}{2}\|D^\ell\bfU(t)\|_\infty+\Sum{|\ell|=0}{1}\|D^\ell{\sf Q}(t)\|_\infty\Big)=0\,.
\eeq{JsB}
From \eqref{def} and the $T$-periodicity of $\bfw_1$ we have for all $n\in\nat$
$$
\bfw_1(x,t)=\bfw_1(x,t+n\,T)=\bfW_1(y,t+n\,T)=\bfU(y,t+n\,T)+\bfv(y,t+n\,T)\,.
$$
Thus, setting
\be
{\sf G}:= [\!]\bfG[\!]_{\infty,2}+\Sum{|k|=0}{1}[\!]D^k\nabla\cdot\bfG[\!]_{\infty,|k|+3}\,,
\eeq{gg}
by \eqref{as} and (\ref{estimunn2}) we get
$$\ba{rl}\medskip
|\bfw_1(x,&\!\!t)|\,(|\bfx|+1)\\&\medskip\le |\bfU(y, t+n\,T)|\,(|\bfy|+M+1)+|\bfv(y,t+n\,T)|\,(|\bfy|+M+1)\\
&\le |\bfU(y, t+n\,T)|\,(|\bfy|+M+1)+C\,{\sf G}\,, 
\ea
$$
and, similarly,
$$\ba{ll}\medskip
|\nabla\bfw_1(x,t)|\,(|\bfx|^2+1)\le |\nabla\bfU(y, t+n\,T)|\,(|\bfy|^2+M+1)+C\,{\sf G}\,.
\\
|D^2\bfw_1(x,t)|\,(|\bfx|^3+1)\le |D^2\bfU(y, t+n\,T)|\,(|\bfy|^3+M+1)+C\,{\sf G}
\ea
$$
Thus, if we pass to the limit $n\to\infty$ in the relations above and use \eqref{A}, \eqref{B}, \eqref{gg},\eqref{bfrt}, and \eqref{JsB} we conclude
\be
\Sum{|k|=0}{2}[\!]D^k\bfw_1[\!]_{\infty,|k|+1}\le C\,\mathscr D_2\,.
\eeq{end}
We now recall that $\bfu=(1-\psi)\bfu+\bfw_1$, and so
the claimed asymptotic property of $\bfu$ follows \eqref{Li_2} and  \eqref{end}. We next observe that from \eqref{sa1}$_1$, \eqref{def}$_2$, \eqref{saz}$_2$ and $T$-periodicity, we get
$$
{\sf p}(x,t)+\tilde{\sf p}(x,t)= {\sf p}(x,t+nT)+\tilde{\sf p}(x,t+nT)={\sf Q}(y,t+nT)+\phi(y,t+nT)\,.
$$
Arguing as in the estimate of $\bfw_1$ and taking into account (\ref{estimunn2}) and \eqref{JsB}, from the preceding relation we deduce
\be
|{\sf p}(x,t)|(|x|^2+1)\le |\tilde{\sf p}(x,t)|(|x|^2+1)+C\,\left([\!]\bfG[\!]_{\infty,2}+[\!]\nabla\cdot\bfG[\!]_{\infty,3}\right)\,.
\eeq{sa2}
Recalling that $h=\nabla\psi\cdot\bfu$ (see \eqref{frt0}$_2$), we infer $\int_{\real^3}\partial_th=0$. Therefore, from \eqref{sa1}$_1$, also after integrating by parts, we deduce
\be\ba{rl}\medskip
\tilde{\sf p}(x,t)&=\Int{\real^3}{}[\cale (x-y)-\cale (x)]\partial_th(y,t)dy\\ &+\Int{\real^3}{}\nabla\cale (x-y)\cdot[\bfxi(t)h(y,t)+\nabla h(y,t)]dy:={\sf I}_1+{\sf I}_2\,.
\ea\eeq{al}
By the mean value theorem and \eqref{GE}, we infer
\be
|{\sf I}_1|\le c_1\,|x|^{-2}\,\|\partial_th(t)\|_1\,,\ \ |{\sf I}_2|\le c_2(\xi_0)\,|x|^{-2}\,\|h(t)\|_{1,1}\,,\ \ \,|x|\ge 4\rho\,,
\eeq{al1}
whereas from classical results on convolutions (e.g. \cite[Theorem II.11.2]{GaB})
\be
|{\sf I}_1|\le c_3\,\|\partial_th(t)\|_q\,,\ \ |{\sf I}_2|\le c_4(\xi_0)\,\|h(t)\|_{1,q}\,,\ \ \,|x|\le 4\rho\,,\ \ q>3\,.
\eeq{al2}
Since
$$
\|\partial_th(t)\|_{L^\infty(L^6)}+\|h(t)\|_{L^\infty(W^{1,6})}\le c\,\big(\|\partial_t\bfu(t)\|_{L^\infty(L^6)}+\|\bfu(t)\|_{L^\infty(W^{1,6})}\big)
$$
from the latter,  \eqref{sa2}--\eqref{al2}, classical embedding,  \eqref{estik}, and  \eqref{A}, \eqref{bfrt}, we conclude
\be
[\!]{\sf p}[\!]_{\infty,2}\le C\,\mathscr D_2\,.
\eeq{al3}
Now, as before, we recall that $p=(1-\psi)p+{\sf p}$, so that from \eqref{Li_1}    and \eqref{al3} we prove the desired property for $p$. In an entirely analogous way one can deduce the pointwise estimate for $\nabla p$.
The proof of the proposition is therefore completed.\QED

\setcounter{equation}{0}
\section{Unique Solvability of the Nonlinear Problem}
We introduce  the following function class:
$$\ba{ll}\medskip
\mathscr X:=\big\{\mbox{$T$-periodic $\bfu$}: \bfu\in  W^{2,2}(D^{2,2})\cap W^{1,2}(D^{4,2})\cap W^{2,\infty}(W^{1,2})\cap L^\infty(D^{3,2})\,;\\ \hspace*{3.5cm} \sum_{|k|=0}^2 [\!]D^k\bfu[\!]_{\infty,|k|+1}<\infty\,;\ \Div\bfu=0\big\}\,.  
\ea$$
Clearly, $\mathscr X$ becomes a Banach space when endowed with the norm
\be
\|\bfu\|_{\mathscr X}:=\|\bfu\|_{W^{2,2}(D^{2,2})}+\|\bfu\|_{W^{1,2}(D^{4,2})}+\|\bfu\|_{W^{2,\infty}(W^{1,2})}+\|\bfu\|_{L^{\infty}(D^{3,2})}
+\sum_{|k|=0}^2[\!]\bfu[\!]_{\infty,|k|+1}\,.
\eeq{00}
Moreover, we set
$$
\mathscr P:=\{\mbox{$T$-periodic}\ p:\ p\in L^\infty(W^{1,2})\cap W^{1,2}(D^{3,2})\,;\ \sum_{|k|=0}^1[\!]D^k p[\!]_{\infty,|k|+1}<\infty\big\}\,,
$$
with
$$
\|p\|_{\mathscr P}:= \|p\|_{L^\infty(W^{1,2})}+\|p\|_{W^{1,2}(D^{3,2})}+\sum_{|k|=0}^1[\!]D^k p[\!]_{\infty,|k|+1}\,.
$$
The main result of this section reads as follows.
\Bt Let $\bfxi\in W^{3,2}(0,T)$ be $T$-periodic with $\int_0^T\bfxi(t)dt=\0$. Moreover, suppose that $\bfb=\Div\bfB$, where $\bfB$ is a $T$-periodic tensor function such that
$$
\bfb\in  W^{2,2}(L^{2})\cap W^{1,2}(D^{2,2})\,,\ 
\|\nabla\bfB\|_{L^\infty(L^\infty(\Omega_{2\rho})}+[\!]{\bfB}[\!]_{\infty,2}+\Sum{|k|=0}{1} [\!]D^k\bfb[\!]_{\infty,|k|+3}<\infty\,,
$$
for some fixed $\rho>R_*$.
Then, setting 
\be\ba{rl}\smallskip
\mathtt D:=\|\bfb\|_{W^{2,2}(L^2)}+\|\partial_t\bfb\|_{L^2(D^{2,2})}&+\Sum{|k|=0}{1} [\!]D^k\bfb[\!]_{\infty,|k|+3}\\&+\|\nabla\bfB\|_{L^\infty(L^\infty(\Omega_{2\rho})}+[\!]\bfB[\!]_{\infty,2}+\|\bfxi\|_{W^{3,2}(0,T)}\,,\ea
\eeq{Dtt}
there exists $\varepsilon_0>0$ such that if 
$\mathtt D<\varepsilon_0\,,
$
problem \eqref{PR} has one and only one  solution $(\bfu,p)\in \mathscr X\times \mathscr P$. Moreover, this solution obeys the following inequality
\be
\|\bfu\|_{\mathscr X}+\|p\|_{\mathscr P}\le C\, \mathtt D\,.
\eeq{ineq}
\ET{m}
{\em Proof.} We employ the contraction mapping theorem. To this end, define the map
$$
M:\textbf{\textsf{u}}\in\mathscr X\mapsto \bfu\in\mathscr X\,,
$$
with $\bfu$ solving the linear problem
\be\ba{cc}\smallskip\left.\ba{ll}\medskip
{\partial}_t\bfu-\bfxi(t)\cdot\nabla\bfu=\Delta\bfu-\nabla {p}+\textbf{\textsf{u}}\cdot\nabla \textbf{\textsf{u}}+\bfb\\
\Div\bfu=0\ea\right\}\ \ \mbox{in $\Omega\times (0,T)$}\\
\bfu(x,t)=\bfxi(t)\,,\ \ (x,t)\in \partial\Omega\times [0,T]\,,
\ea
\eeq{lin}
Set 
\be\bfsff:=\textbf{\textsf{u}}\cdot\nabla\textbf{\textsf{u}}=\Div(\textbf{\textsf{u}}\otimes \textbf{\textsf{u}}):=\Div\bfsfF\,,
\eeq{C}
where we used the condition $\Div\textbf{\textsf{u}}=0$. Clearly,
\be\ba{ll}\medskip
\|\bfsff\|_{W^{2,2}(L^2)}\le c\,\Big(\|\textbf{\textsf{u}}\cdot\nabla\textbf{\textsf{u}}\|_{L^2(L^2)}+\|\partial_t\textbf{\textsf{u}}\cdot\nabla\textbf{\textsf{u}}\|_{L^2(L^2)}+\|\textbf{\textsf{u}}\cdot\nabla\partial_t\textbf{\textsf{u}}\|_{L^2(L^2)}\\ \medskip
\hspace*{3cm}+\|\partial_t^2\textbf{\textsf{u}}\cdot\nabla\textbf{\textsf{u}}\|_{L^2(L^2)}+\|\partial_t\textbf{\textsf{u}}\cdot\nabla\partial_t\textbf{\textsf{u}}\|_{L^2(L^2)}+\|\textbf{\textsf{u}}\cdot\nabla\partial_t^2\textbf{\textsf{u}}\|_{L^2(L^2)}\Big)\\ \medskip
\|\bfsff\|_{W^{1,2}(D^{2,2})}\le c\, \Big(
\||D^2\textbf{\textsf{u}}|\,|\nabla \textbf{\textsf{u}}|\|_{L^2(L^2)}^2+\|\textbf{\textsf{u}}\cdot\nabla D^2\textbf{\textsf{u}}\|_{L^2(L^2)}\\
\hspace*{4cm}+\|D^2\big(\partial_t\textbf{\textsf{u}}\cdot\nabla\textbf{\textsf{u}}+\textbf{\textsf{u}}\cdot\nabla\partial_t\textbf{\textsf{u}}\big)\|_{L^2(L^2)}\Big)
\,.
\ea
\eeq{D}
Thus, by a straightforward calculation, we show
\be\ba{rl}\medskip
\|\bfsff\|_{W^{2,2}(L^2)}\le&\! c\,\big[[\!]\bfsfu[\!]_{\infty,1}\big([\!]\nabla\bfsfu[\!]_{\infty,2}+\|\bfsfu\|_{W^{2,2}(W^{1,2})}\big)\\&\!+[\!]\nabla\bfsfu[\!]_{\infty,2}\|\bfsfu\|_{W^{2,2}(L^2)}+\|\partial_t\bfsfu\|_{L^\infty(L^\infty)}\|\bfsfu\|_{W^{1,2}(W^{1,2})}\big]\,.\ea 
\eeq{E}
Employing in the last term of \eqref{E} the classical embedding inequality:
\be
\|\partial_t\bfsfu\|_{L^\infty(L^\infty)}\le c\,\big(\|\partial_t\bfsfu\|_{L^\infty(L^2)}+\|\partial_t\bfsfu\|_{W^{1,2}(D^{2,2})}\big)\,,
\eeq{ei}
from \eqref{E} and \eqref{00} we then conclude
\be
\|\bfsff\|_{W^{2,2}(L^2)}\le c\,\|\bfsfu\|_{\mathscr X}^2\,.
\eeq{vb1}
In a similar fashion, we show
\be\ba{rl}\medskip
\|\bfsff\|_{W^{1,2}(D^{2,2})}\le &\!c\,\big[[\!]D^2\bfsfu[\!]_{\infty,3}\big([\!]\nabla\bfsfu[\!]_{\infty,2}+\|\bfsfu\|_{W^{1,2}(W^{1,2})}\big)\\
&+\|\bfsfu\|_{W^{1,2}(D^{3,2})}\big([\!]\bfsfu[\!]_{\infty,1}+\|\partial_t\bfsfu\|_{L^\infty(L^\infty)}\big)+[\!]\nabla\bfsfu[\!]_{\infty,2}\|\bfsfu\|_{W^{1,2}(D^{2,2})}\big]\,.
\ea\eeq{F}
Again, by classical embedding,
\be
\|\bfsfu\|_{W^{1,2}(D^{3,2})}\le c\,\big(\|\bfsfu\|_{W^{1,2}(D^{2,2})}+\|\bfsfu\|_{W^{1,2}(D^{4,2})}\big)\,.
\eeq{fesso}
Therefore, employing \eqref{ei} and \eqref{fesso} in \eqref{F} and taking into account \eqref{00} we deduce
\be
\|\bfsff\|_{W^{1,2}(D^{2,2})}\le c\,\|\bfsfu\|_{\mathscr X}^2\,.
\eeq{fico}
Finally, and obviously,
\be\ba{rl}\medskip
[\!]\bfsfF[\!]_{\infty, 2}+\Sum{|k|=0}1[\!]D^k\bfsff[\!]_{\infty, |k|+3}+\|\nabla\bfsfF\|_{L^\infty(L^\infty(\Omega_{2\rho})}&\!\le c\,\big([\!]\textbf{\textsf{u}}[\!]_{\infty,1}\Sum{|k|=0}{2}[\!]D^k\bfsfu[\!]_{|k|+1}+[\!]\nabla\bfsfu[\!]_{\infty,2}^2\big)\\
&\! \le c\,\|\bfsfu\|_{\mathscr X}^2\,.\ea
\eeq{G}
As a result, from \eqref{C}, \eqref{vb1}, \eqref{fico} and \eqref{G} we find that  $\bfsff$ and $\bfsfF$
satisfy the assumptions of $\bff$ and $\bfcalf$, respectively, in \propref{1}, and, in addition,
\be
\|\bfsff\|_{W^{2,2}(L^2)}+\|\bfsff\|_{W^{1,2}(D^{2,2})}+\|\bfsfF\|_{2,\Omega_{2\rho}}+\Sum{|k|=0}1[\!]D^k\bfsff[\!]_{\infty,|k|+3}\le c\, \|\bfsfu\|_{\mathscr X}^2\,.
\eeq{H}
Thus, by that proposition and the assumption on the data, we deduce, on the one hand, that $\big(M(\textbf{\textsf{u}}),p)\in\mathscr X\times\mathscr P$ --so that, in particular, $M$ is well defined-- and, on the other hand, that $\bfu=M(\textbf{\textsf{u}})$ obeys the estimate:
\be
\|\bfu\|_{\mathscr X}+\|p\|_{\mathscr P}\le C_1\,\big(\|\bfu\|_{\mathscr X}^2+\mathtt D\big)\,.
\eeq{I}
Next, suppose $\|\textbf{\textsf{u}}\|_{\mathscr X}<\delta$. From \eqref{I}  it  follows
$$
\|\bfu\|_{\mathscr X}\le C_1\,\big(\delta^2+\mathtt D\big)\,,
$$
from which we infer that if we pick
\be
\mathtt D<\min\{\frac{\delta}{2 C_1},\xi_0\}\,, \ \,\delta<\frac1{2C_1}\,,
\eeq{L}
we obtain
\be
\|\bfu\|_{\mathscr X}< \delta\,.
\eeq{M}
Let $\textbf{\textsf{u}}_i\in\mathscr X$ $i=1,2$, and set 
$$
\textbf{\textsf{u}}:=\textbf{\textsf{u}}_1-\textbf{\textsf{u}}_2\,,\ \ \bfu:=M(\textbf{\textsf{u}}_1)-M(\textbf{\textsf{u}}_2)\,.
$$
From \eqref{lin} we then get
\be\ba{cc}\smallskip\left.\ba{ll}\medskip
{\partial}_t\bfu-\bfxi(t)\cdot\nabla\bfu=\Delta\bfu-\nabla {p}+\textbf{\textsf{u}}_1\cdot\nabla \textbf{\textsf{u}}+\textbf{\textsf{u}}\cdot\nabla\textbf{\textsf{u}}_2\\
\Div\bfu=0\ea\right\}\ \ \mbox{in $\Omega\times (0,T)$}\\
\bfu(x,t)=\0\,,\ \ (x,t)\in \partial\Omega\times [0,T]\,.
\ea
\eeq{line}
Arguing as in the proof of \eqref{H} we can show
$$
\|\bfu\|_{\mathscr H}\le C_1\,\left(\|\textbf{\textsf{u}}_1\|_{\mathscr X}+\|\textbf{\textsf{u}}_2\|_{\mathscr X}\right)\|\textbf{\textsf{u}}\|_{\mathscr X}\,.
$$
Consequently, if $\|\textbf{\textsf{u}}_i\|_{\mathscr X}<\delta$, $i=1,2$, from the preceding inequality we find
$$
\|\bfu\|_{\mathscr H}< 2C_1\delta\|\textbf{\textsf{u}}\|_{\mathscr X}\,,
$$
and since by \eqref{L} $2C_1\delta<1$, we may conclude that $M$ is a contraction, which ends the proof of existence. Finally, the estimate \eqref{ineq} is a consequence of \eqref{I}, \eqref{M} and the choice of $\delta$ in \eqref{L}.\QED 
\begin{remark} 
%
In the particular case $\bfxi(t)\equiv\0$, \theoref{m} furnishes (in a better regularity class and with more information about the behavior at infinity)  existence results similar to those proved \cite{GaSo,KMT}.
\end{remark}
\Br
\theoref{m} establishes the uniqueness of the solution in the ball of $\mathscr X$ of radius $\delta$. However, a more general uniqueness result ``in the large" could be actually shown in a sufficiently regular class of solutions (not necessarily ``small"), and even in a suitable class of ``weak" solutions. In fact, the former could be attained by employing the same ``cut-off" procedure used in the proof of \lemmref{4.2}, in conjunction with the pointwise asymptotic properties of the solution constructed in \theoref{m}. As for the latter, one could just follow, step by  step,  the proof provided in \cite[Theorem 5]{GaSo}.     
\ER{2}\setcounter{equation}{0}
\section{Asymptotic Spatial Behavior and Steady Streaming} 
\theoref{m} asserts, in particular, that $\bfu$,  $p$ and some of their derivatives have a polynomial (spatial) decay rate at large distance from the body $\mathscr B$. Objective of this section is to provide a more detailed analysis of this property and show that, ``far" from $\mathscr B$, the flow velocity field presents  a distinctive steady-state character, in spite of being driven by a time-periodic mechanism. This rigorous finding is in agreement with the classical  phenomenon of ``steady streaming"  observed in the motion of a viscous liquid past an oscillating body; see \cite[p. 428--432]{Sch}, \cite{Ri} and the references therein. 
\smallskip\par
To prove the above, we recall the following splitting of $\bfu$ into its averaged and oscillatory components (see  \eqref{471}):
$$
\bfu=\bar{\bfu}+\bfw\,,
.
$$
The following lemma holds.
\Bl The oscillatory component $\bfw$ of the solution $\bfu$ of \theoref{m} satisfies
$$
[\!]\bfw[\!]_{2,\infty}<\infty\,.
$$
If, in particular, $\bfxi(t)\equiv\0$, then the faster decay condition is valid:
$$
[\!]\bfw[\!]_{\infty, 3}<\infty.
$$
\EL{4.1} 
{\em Proof.} Since $\bar{\bfw}=\0$, from the Poincar\'e inequality we get,  for all $x\in\Omega$:
$$
\int_0^T|\bfw(x,t)|dt\le T\int_0^T|\partial_t\bfw(x,t)|dt\,,
$$
which once combined with an elementary embedding inequality, implies
\be
\sup_{s\in [0,T]}|\bfw(x,s)|\le c\int_0^T|\partial_t\bfw(x,t)|dt\,.
\eeq{4.1}
Therefore, the claimed properties follow directly  from \eqref{4.1}, \eqref{PR}$_1$ and the pointwise decay  estimates established in \theoref{m}.\QED

We next observe that, from \eqref{PR}, the averaged component, $\bar{\bfu}$, of $\bfu$ and corresponding averaged pressure $\bar{p}$ solve the following boundary-value problem
\be\ba{cc}\medskip\left.\ba{ll}\medskip
\Delta\bar{\bfu}=\bar{\bfu}\cdot\nabla\bar{\bfu}+\nabla\bar{p}-\Div\bfF\\
\Div\bar{\bfu}=0\ea\right\}\ \ \mbox{in\,$\Omega$}\,,\\
\bar{\bfu}=\0\ \ \mbox{on $\partial\Omega$}\,, 
\ea
\eeq{bvp}
where
\be
\bfF:=\big(\bar{\bfw\otimes\bfw}-\bar{\bfB}\big)-\bar{\bfxi\otimes\bfw}:=
\bfF_1+\bfF_2\,.
\eeq{FU}
\begin{definition} Two vector fields $\bfU_1$, $\bfU_2\in L^\infty(\Omega)$ are {\em asymptotically equivalent} --and we write $\bfU_1\sim \bfU_2$-- if (i) $[\!]\bfU_i[\!]_1<\infty$, $i=1,2$, and (ii) $[\!]\bfU_1-\bfU_2[\!]_{1+\delta}<\infty\,,$ for some $\delta>0$.  
\end{definition}

\Bl Let $\bfG\in L^\infty(\real^3)$ with support in $B_R$. Further, let $\psi=\psi(|x|)$ be a smooth function that is 0 in $B_{R/2}$ and 1 in $B^R$, $R\ge 2R_*$. There exists $\varepsilon_1>0$ such that if
\be
[\!]\bfF_2[\!]_{2}+\|\bfG\|_\infty\le \varepsilon_1\, 
\eeq{chm}
then the problem
\be
\left.\ba{ll}\medskip
\Delta \bfU=\bfU\cdot\nabla\bfU+\nabla P-\Div(\psi\bfF_2)+\bfG\\
\Div\bfU=0\ea\right\}\ \ \mbox{in $\real^3$}\,,
\eeq{MH1}
has at least one solution $(\bfU,P)\in W^{2,2}_{\rm loc}(\real^3)\times W^{2,1}_{\rm loc}(\real^3)$ with $[\!]\bfU[\!]_{1}<\infty$, and, moreover, 
\be
[\!]\bfU[\!]_{1}\le C\,\left([\!]\bfF_2[\!]_{2}+\|\bfG\|_\infty\right)
.
\eeq{MH0} 
Finally, let $(\bfU',P')$ solve the problem 
\be\left.\ba{ll}\medskip
\Delta \bfU'=\bfU'\cdot\nabla\bfU'+\nabla P'-\Div(\psi\bfF_2)+\bfG'\\
\Div\bfU'=0\ea\right\}\ \ \mbox{in $\real^3$}\,,
\eeq{MH2}
with $\bfG'$ satisfying the same properties listed for $\bfG$. Then, if
\be
\int_{B_R}\bfG=\int_{B_R}\bfG',
\eeq{avem}
we have $\bfU\sim\bfU'$.
\EL{4.2_1}
{\em Proof.} We begin to notice that, in view of \lemmref{4.1}, the assumption $[\!]\bfF[\!]_2<\infty$ is meaningful. Let 
$$
\bfH(x)=\int_{\real^3}\nabla\cale(x-y)\cdot\bfG(y)dy\,.
$$
Then, clearly $\Div\bfH=\bfG$. Furthermore, proceeding as in \eqref{ha1}, \eqref{ha2}, we show
$$
[\!]\bfH[\!]_2\le c\,\|\bfG\|_\infty\,,
$$
so that, in view of our assumptions, the field $\bfF+\bfH$ meets the hypotheses of \cite[Lemma X.9.1]{GaB}. As a result, there is a corresponding solution $(\bfU,P)$ to \eqref{MH1} such that
\be
(\bfU,P)\in D_0^{1,q}(\real^3)\times L^q(\real^3)\,,\ \ \mbox{all $q>3/2$}\,,
\eeq{reg} 
which, in addition, satisfies $[\!]\bfU[\!]_1<\infty$, along with the estimate \eqref{MH0}. Finally, from \eqref{reg}   and classical regularity results \cite[Theorem X.1.1]{GaB} we infer $(\bfU,P)\in W^{2,2}_{\rm loc}(\real^3)\times W^{2,1}_{\rm loc}(\real^3)$. Next, setting 
$$
\bfzeta:=\bfU-\bfU'\,,\ \ \bfg:=\bfG-\bfG'\,,
$$
from \eqref{MH1}--\eqref{avem} we deduce the following integral representation \cite[Theorem X.5.2]{GaB}:
\be
\zeta_i(x)=\int_{\real^3}\big[S_{ij}(x-y)-S_{ij}(x)\big]g_j(y)dy+\int_{\real^3}\pde{}{x_\ell}S_{ij}(x-y)\big[\zeta_\ell U_j-U_\ell'\zeta_j\big](y)dy\,,
\eeq{mo}
where $\bfS$ is the Stokes fundamental tensor that, we recall, satisfies the following asymptotic bounds \cite[Section IV.2]{GaB}
\be
|D^k\bfS(\chi)|\le C\,|\chi|^{-(1+|k|)}\,,\ \ \bfchi\neq \0\,,\ \ |k|\in \nat\cup\{0\}\,.
\eeq{mo1}
We now regard \eqref{mo1} as an integral equation in the unknown $\bfzeta$. It is simple to show that  this equation has a solution, $\hat{\bfzeta}$, in the space
$$
\mathscr S_\alpha:=\{\bfz\in L^\infty(\real^3):\ [\!]\bfz[\!]_{1+\alpha}<\infty\}\,, 
$$
for some $\alpha\in(0,1)$, provided we take $\varepsilon_1$ appropriately ``small." Actually, recalling that $\supp(\bfg)\subset B_R$ and that both $\bfU$, $\bfU'$ are in $\mathscr S_0$, from \eqref{mo} and \eqref{mo1} we show
\be
|\bfzeta(x)|\le C_1\,\|\bfg\|_\infty(1+|x|)^{-2}+C_2[\!]\bfzeta[\!]_{1+\alpha}\big([\!]\bfU[\!]_{1}+[\!]\bfU'[\!]_{1}\big)\int_{\real^3}\frac{dy}{|x-y|^2|y|^{2+\alpha}}\,.
\eeq{mo2}
On the other hand, by  \cite[Lemmas II.9.2, II.11.2]{GaB} we have
\be
\int_{\real^3}\frac{dy}{|x-y|^2|y|^{2+\alpha}}\le C(\alpha)(1+|x|)^{-1-\alpha}\,,\ \ \alpha\in [0,1)\,.
\eeq{mo3}
Thus, using \eqref{mo2} and \eqref{mo3},  by a simple contraction argument it follows that, for a given $\alpha\in (0,1)$, we can choose a corresponding $\varepsilon_1$ in \eqref{chm} such that \eqref{mo} has a solution $\hat{\bfzeta}\in \mathscr S_\alpha$. It is also readily proved that $\hat{\bfzeta}=\bfzeta$. In fact, setting $\bfz:=\bfzeta-\hat{\bfz}$, we have
$$
z_i(x)=\int_{\real^3}\pde{}{x_\ell}S_{ij}(x-y)\big[z_\ell U_j-U_\ell'z_j\big](y)dy\,,
$$
and so, employing in this relation \eqref{mo1}, \eqref{mo3} with $\alpha=0$, \eqref{chm}, and \eqref{MH0} we get
$$
[\!]\bfz[\!]_1\le C\,\varepsilon_1 [\!]\bfz[\!]_1\,,
$$
which allows us to conclude $\bfz\equiv\0$ by taking $\varepsilon_1$ sufficiently small. The proof is completed.\QED
\Br
If $\bfxi\equiv\0$, namely,  $\bfF_2\equiv \0$, then every solution in \lemmref{4.2_1} corresponding to some $\bfG$ that satisfies the assumption of that lemma along with the condition
$$
\int_{B_R}\bfG:=\bfbeta\neq\0\,,
$$
is asymptotically equivalent to  a specific member of the well-known Landau family of solution \cite{L,KS}. To see this, let  $(r,\theta,\phi)$ be a system of polar coordinates, with polar axis oriented in the direction $\bfbeta/|\bfbeta|$ which, without loss, we take coinciding with the positive $x_1-$direction. We recall that the {\em Landau solution}   corresponding to $\bfbeta$ is a pair $(\bfU^{\beta},P^{\beta})$ satisfying
\be\left.\ba{ll}\medskip
\Delta \bfU^\beta-\bfU^\beta\cdot\nabla\bfU^\beta-\nabla P^\beta=\bfbeta\,\delta\\
\Div\bfU^\beta=0\ea\right\}\ \ \mbox{in $\real^3$}\,,
\eeq{LS}
with $\delta$ Dirac distribution, and
defined,  for $r>0$, as follows
\be\ba{rl}\medskip
U^{\beta}_r &= \Frac{2}{r}\left[\Frac{A^2-1}{(A-\cos\theta)^2}-1\right]\,,\\ \medskip
U^{\beta}_\theta &= - \Frac{2 \sin \theta}{r(A-\cos\theta)}\,,\\ \medskip
U^{\beta}_\phi&=0\,,\\ 
P^{\beta} &= \Frac{ 4(A\cos\theta-1)}{
r^2(A-\cos \theta)^2} \,,
\ea
\eeq{IXv.9.LN}
where the parameter $A\in(1,\infty)$ is chosen in such a way that
\be
16\pi\left(A +
\frac{1}{2}A^2 \log\frac{A -1}{A + 1}
+\frac{4A}{3(A^2- 1)}
\right)=\beta
\eeq{IXv.9sbN}
Since the function on the left-hand side is monotonically decreasing in $A\in (1,\infty)$ and its range coincides with $(0,\infty)$,  we deduce that for any given $\beta$ ($>0$) there is one and only one $A$ satisfying \eqref{IXv.9sbN}, namely, one and only one Landau solution $(\bfU^{\beta},P^{\beta})$. Moreover, observing that $A\to\infty$ as $\beta\to 0$, from \eqref{IXv.9.LN} we also deduce, in particular, 
\be
[\!]\bfU^\beta[\!]_1+
[\!]\nabla\bfU^\beta[\!]_2+
[\!]P^\beta[\!]_2\le\kappa (\beta)\,,\ \ \lim_{b\to0}\kappa(\beta)=0\,.
\eeq{MH}
Now, by following a standard procedure, we regularize $(\bfU^\beta,P^\beta)$ around $x=0$ by defining $\tilde{\bfU^\beta}:=\psi\bfU^\beta-\textbf{\textsf{U}}$, $\tilde{P^\beta}:=\psi\,P^\beta$, where $\psi=\psi(|x|)$ is the ``cut-off"  function introduced in \lemmref{4.2_1}, while $\Div \textbf{\textsf{U}}=\nabla\psi\cdot\bfU^\beta$ in $B_R$, $\textbf{\textsf{U}}\in C_0^\infty(B_R))$ \cite{KS}. It is then readily checked that $(\tilde{\bfU^\beta},\tilde{P^\beta})$ is a solution to the following problem
\be
\left.\ba{ll}\medskip
\Delta \tilde{\bfU^\beta}=\tilde{\bfU^\beta}\cdot\nabla\tilde{\bfU^\beta}+\nabla \tilde{P^\beta}+\tilde{\bfG^\beta}\\
\Div\tilde{\bfU^\beta}=0\ea\right\}\ \ \mbox{in $\real^3$}\,,
\eeq{MH_0}
with $\tilde{\bfG^\beta}\in C_0^\infty(B_R)$ and such that, by \eqref{LS} and \eqref{MH},
\be
\int_{B_R}\tilde{\bfG^\beta}=\beta\,\bfe_1\,,\ \ \|\tilde{\bfG^\beta}\|_\infty\le C\,\kappa(\beta)\;
\eeq{MH_00}
see \cite{KS} for details. Thus, the claimed asymptotic equivalence is a consequence of \lemmref{4.2}.\label{rem:landau}
\ER{landau}
\par
We recall the definition of Cauchy stress:
$$
\bfT(\textbf{\textsf{w}},{\sf p}):=-{\sf p}\,\bfI+\nabla \textbf{\textsf{w}}+(\nabla \textbf{\textsf{w}})^\top\,,
$$
with $\bfI$ identity matrix and $\top$ denoting transpose.
\Bl Let $(\bar{\bfu},\bar{p})$ be the averaged component of the solution $(\bfu,p)$ given in \theoref{m},  satisfying \eqref{bvp}, \eqref{FU}, and let $(\bfU,P)$ be the solution to \eqref{MH1} constructed in \lemmref{4.2_1},  corresponding to a vector field 
\be
\bfG=\theta(x)\int_{\partial\Omega}\big(\bfT(\bar{\bfu},\bar{p})-\bar{\bfB}\big)\cdot\bfn\,,\ \ \theta\in C_0(B_R)\,,\ \ \int_{B_R}\theta=1\,. 
\eeq{mitt}
Then, there is $\varepsilon_2>0$ such that, if 
$${\mathtt D}+[\!]\bar{\bfB}[\!]_{2+\alpha}<\varepsilon_2
$$ 
for some $\alpha\in (0,1)$  {\rm(}$\mathtt D$ defined in \eqref{Dtt}{\rm)},  necessarily $\bar{\bfu}\sim \bfU$.
\EL{4.3}
{\em Proof.} Let $\psi=\psi(|x|)$ be the ``cut-off" function introduced in \lemmref{4.2_1}, and set
$$
\bfv:=\psi\,\bar{\bfu}-\textsf{\textbf{V}}\,,\ \ q:=\psi\,p\,,
$$
where $\bfV\in W_0^{4,2}(B_R)$ satisfies $\Div\textsf{\textbf{V}}=\nabla\psi\cdot\bar{\bfu}$ with 
\be
\|\textsf{\textbf{V}}\|_{4,2}\le c\,\|\bar{\bfu}\|_{3,2,\Omega_R}\,;
\eeq{SBV0}
see \cite[Theorem III.3.3]{GaB}. In view of  \eqref{bvp}, \eqref{FU} and  the regularity properties of $\bfu$, we show that $(\bfv,q)$ satisfies the following problem
\be
\left.\ba{ll}\medskip
\Delta \bfv=\bfv\cdot\nabla\bfv+\nabla q-\Div(\psi\bfF_1)-\Div(\psi\bfF_2)+\bfG_0\\
\Div\bfv=0\ea\right\}\ \ \mbox{in $\real^3$}\,,
\eeq{mono} 
where $\bfG_0\in L^\infty(B_R)$ with $\supp(\bfG_0)\subset B_R$. Moreover, also using \eqref{SBV0}, we have
\be
\|\bfG_0\|_\infty\le C\,\big(\|\bar{\bfu}\|_{3,2,\Omega_R}^2+\|\bar{\bfu}\|_{3,2}+\|\bar{\bfu}\|_{1,\infty}+\|\bar{p}\|_{\infty}+
\|\bar{(\bfw-\bfxi)\otimes\bfw}\|_\infty+\|\bar{\bfB}\|_\infty\big)\,.
\eeq{SBV}
Observing that
$$
\Div\big[\bfT(\bfv,q)-\bfv\otimes\bfv+\psi(\bfF_1+\bfF_2)\big]= \Delta \bfv-\bfv\cdot\nabla\bfv-\nabla q+\Div(\psi\bfF_1)+\Div(\psi\bfF_2)\,,
$$
and that $\psi(R)=1$, integrating both sides of \eqref{mono}$_1$ over $B_R$, we get
\be\ba{rl}\medskip
\Int{B_R}{}\bfG_0=&\!\Int{\partial B_R}{}\big[\bfT(\bar{\bfu},\bar{p})-\bar{\bfu}\otimes\bar{\bfu}+\bfF_1+\bfF_2\big]\cdot\bfn\\
=&\Int{\partial \Omega}{}\big[\bfT(\bar{\bfu},\bar{p})-\bar{\bfB}\big]\cdot\bfn\,,
\ea\eeq{sfax}
where we have used the following properties, consequences of \eqref{bvp}$_{1,3}$:
$$\ba{ll}\medskip
\Div\big[\bfT(\bar{\bfu},\bar{p})-\bar{\bfu}\otimes\bar{\bfu}+ \bfF_1+\bfF_2\big]=\0\\
\Int{\partial\Omega}{}\big(\bar{\bfw\otimes \bfw-\bfxi\otimes\bfw}\big)\cdot\bfn=\Int{\partial\Omega}{}\big(\bar{\bfxi\otimes \bfxi-\bfxi\otimes\bfxi}\big)\cdot\bfn=\0\,.
\ea
$$
Therefore,
$$
\int_{B_R}\bfG=\int_{B_R}\bfG_0\,.
$$
Let $\bfzeta:=\bfv-\bfU$. starting with  \eqref{mono}, \eqref{MH1} and proceeding as in the proof of \lemmref{4.2_1} (see \eqref{mo}),
we show 
\be
\zeta_i(x)=\int_{\real^3}\big[S_{ij}(x-y)-S_{ij}(x)\big]g_j(y)dy+\int_{\real^3}\pde{}{x_\ell}S_{ij}(x-y)\big[\zeta_\ell v_j-U_\ell\zeta_j-\psi\,F_{1\ell j}\big](y)dy\,,
\eeq{mono_1}
where $\bfg:=\bfG_0-\bfG$. Since,  by  classical trace theorems  and \lemmref{4.1}, one shows
$$\ba{rl}\medskip
\|\bfG\|_\infty+&\![\!]\bfF_2[\!]_{2}+[\!]\bfF_1[\!]_{2+\alpha}\\  &\le C_1\big(\|\bar{\bfu}\|_{2,2,\Omega_R}+\|\bar{p}\|_{1,2,\Omega_R}+[\!]\bfw[\!]_{\infty,2}+[\!]\bfw[\!]_{\infty,2+\alpha}^2+[\!]\bar{\bfB}[\!]_{2+\alpha}\big)\,,
\ea
$$
from the latter, \eqref{SBV} and \eqref{ineq} we deduce 
$$
\|\bfG_0\|_\infty+\|\bfG\|_\infty+[\!]\bfF_2[\!]_{2}+[\!]\bfF_1[\!]_{2+\alpha}\le C_2\,(\mathtt D+\mathtt D^2)<C\,\varepsilon_2\,.
$$
Thus, we can argue exactly as in the proof of \lemmref{4.2_1} to show $[\!]\bfzeta[\!]_{1+\alpha}<\infty$, which completes the proof of the lemma.\QED
\par
We are now in a position to show the main result of this section.
\Bt Let $(\bfu\equiv \overline{\bfu}+\bfw,p)$ be the solution determined in \theoref{m}. Then, under the assumptions on $\mathtt D$ and $\bar{\bfB}$ of \lemmref{4.3}, the velocity field $\bfu$ has the following representation
$$
\bfu(x,t)=\bfU(x)+\bfsigma(x)+\bfw(x,t)
$$
where $(\bfU,P)$ is the steady-state solution of \lemmref{4.3}, and for some $\alpha\in (0,1)$,
$$
[\!]\bfsigma[\!]_{1+\alpha}+[\!]\bfw[\!]_{\infty,2}<\infty\,.
$$
The field $\bfU$ is unique up to an asymptotically equivalent velocity field.
However, if $\bfxi\equiv\0$, then $\bfU$ is uniquely determined and coincides with the Landau solution $(\bfU^\beta,P^\beta)$ where
$$
\bfbeta=\int_{\partial\Omega}\big(\bfT(\bar{\bfu},\bar{p})-\bar{\bfB}\big)\cdot\bfn\,.
$$
Moreover, in such a case, we have a faster decay of the oscillatory component, namely,
$$
[\!]\bfw[\!]_{\infty,3}<\infty\,.
$$
\ET{m1}
{\em Proof.} It is enough to observe that
$$
\bfu=\bfU+(\bar{\bfu}-\bfU)+\bfw:=\bfU+\bfsigma+\bfw\,,
$$
and employ \lemmref{4.2_1},  \lemmref{4.3} and Remark \ref{rem:landau}.\QED
\begin{remark} In the case $\bfxi\equiv\0$, \theoref{m1} sharpens an analogous result showed in \cite[Theorem 1.2]{KMT}. 

\end{remark}
\textbf{Acknowledgement}.
{Work  partially supported by NSF grant DMS-1614011.}

\ed
set
\be
\mathscr Q_0:\bfv\in\mathscr W_{2\pi,0}^2(\Omega)\mapsto \omega_0\bfv_\tau -{\rm P}\,[\lambda_0\,\partial_1\bfv+\bfe_1\times\bfx\cdot\nabla\bfv-\bfe_1\times\bfv]-\Delta\bfv\in \mathscr H_{2\pi,0}(\Omega)\,.
\eeq{4_2}
\bibitem{MaSo}  Maremonti, P., and Solonnikov, V. A., On nonstationary Stokes problem in exterior domains. {\em Ann. Scuola Norm. Sup. Pisa} Cl. Sci. (4)  {\bf 24}   395--449 (1997)
it remains to show the validity of \eqref{est}$_2$.  To this end, let $h\in L^{q'}(\Omega_R)$, $q'\in[\frac65,\frac32)$, with $\int_{\Omega_R}h=\0$, and let $\bfphi\in W^{1,q'}_0(\Omega_R)$ be a solution to the problem $\Div\bfphi=h$ in $\Omega_R$, satisfying $\|\bfphi\|_{1,q'}\le c_R\|h\|_{q'}$. The existence of such a $\bfphi$ is well known \cite[Theorem III.3.1]{GaB}. Dot-multiplying both sides of \eqref{4.7}$_1$ by $\bfphi$ and integrating by parts over $\Omega_R$, we get
$$
\langle \partial_t\bfu-\bfxi(t)\cdot\nabla\bfu,\bfphi\rangle-\langle\nabla\bfu-\bfcalf,\nabla\bfphi\rangle=\langle p,\Div\bfphi\rangle=\langle p,h\rangle\,.
$$ 
From this relation,  the properties of $\bfphi$ and the arbitrariness of $h$, we deduce that $p$, modified by a possible addition of a ($T$-periodic) function of time, must obeys the following inequality ($q=q'/(q'-1)$)
\be
\|p\|_{q,\Omega_R}\le c\,\big(\|\partial_t\bfu\|_{q,\Omega_R}+\|\nabla\bfu\|_{q,\Omega_R}+\|\bfcalf\|_{q,\Omega_R}\big)\,,\ \ q\in(3,6]\,,
\eeq{Gy}
for all $R\ge R_*$ and with $c=c(R,\xi_0)$; see also \cite[Lemma IV.1.1]{GaB}. By \eqref{class}, we infer, in particular
$$
\partial_t\bfu\in L^\infty(W^{1,2})\,,\
 \ \nabla\bfu\in W^{1,2}(W^{1,2})
$$
so that
from classical  results  it follows $$\partial_t\bfu\in L^\infty(L^6)\,; \ \ \nabla\bfu\in L^\infty(L^{q})\,,\ \ \mbox{ for {\em all} $q\in[2,6]$, with continuous embedding}\,, 
$$
which  with the help of \eqref{est}$_1$ allows us to deduce
\be 
\|\partial_t\bfu\|_{L^\infty(L^6)}+\|\nabla\bfu\|_{L^\infty(L^q)}\le C\,\big(\|\bff\|_{W^{1,2}(L^2)}+\|\bfcalf\|_{L^{2}(L^{2})}+\|\bfxi\|_{W^{2,2}(0,T)} \big)\,,\ \ \mbox{all $q\in [2,6]$}\,.
\eeq{Xm}
The bound \eqref{est}$_2$ then follows from \eqref{Gy}, \eqref{Xm} and the assumption on $\bfcalf$. The proof of existence is thus accomplished.
We then show that if $|\bfxi(t)|$, $|\bfxi'(t)|$ and $|\bfxi''(t)|$, are below a given constant, the corresponding Navier-Stokes problem  has one and only one strong time-periodic solution of period $T$. Moreover, the associated velocity field $\bfv=\bfv(x,t)$  decays like $|x|^{-1}$, while $\nabla\bfv$ and pressure field $p$ decay like $|x|^{-2}$, uniformly in $t\in [0,T]$.

In a similar fashion, by differentiating two times both sides of \eqref{4.8}$_1$ with respect to time and dot-multiplying the resulting equation one time by $\partial^2_t\bfv_k$, a second time by $P\Delta\partial^2_t\bfv_k$, and integrating over $\Omega$, we show
\be
\half\ode{}t\|\partial_t\bfv_k(t)\|_2^2+\|\nabla\partial_t\bfv_k(t)\|_2^2=\langle\bfxi'\cdot\nabla\bfv_k,\partial_t\bfv_k\rangle+\langle\partial_t\bff+\partial_t\bff_c,\partial_t\bfv_k(t)\rangle\,,
\eeq{tem}
and
\be
\ba{rl}\medskip
\half\ode{}t\|\nabla\partial_t^2\bfv_k(t)\|_2^2+\|P&\!\!\Delta\partial_t\bfv_k(t)\|_2^2\\&=\langle\bfxi'\cdot\nabla\bfv_k,P\Delta\partial_t\bfv_k(t)\rangle+\langle\partial_t\bff+\partial_t\bff_c,P\Delta\partial_t\bfv_k(t)\rangle\,.
\ea
\eeq{muo}
\begin{thebibliography}{99}
\bibitem{Gnwp} Galdi, G.P., A steady-state exterior Navier--Stokes problem that is not well-posed. {\em Proc. Amer. Math. Soc.}  {\bf 137}   679--684 (2009) 
\bibitem{GaB}Galdi, G.P., {\em An introduction to the mathematical theory of the Navier-Stokes equations.
Steady-state problems}, Second edition. Springer Monographs in Mathematics,
Springer, New York (2011)
\bibitem{GaA}Galdi, G.P., On time-periodic flow of a viscous liquid past a moving cylinder, {\em Arch. Ration. Mech. Anal.},  {\bf 210} 451--498 (2013)
\bibitem{GaKy}Galdi, G.P., and Kyed, M., Time periodic solutions to the Navier-Stokes
equations, in {\em Handbook of Mathematical Analysis in Mechanics of Viscous Fluids}, Eds Y.Giga, and A. Novotn\'y, Springer-Verlag (2017)
\bibitem{GaKy1} Galdi, G.P., Kyed, M., Time-periodic flow of a viscous liquid past a body.  Partial differential equations in fluid mechanics, 20–49, London Math. Soc. Lecture Note Ser., {\bf 452}, Cambridge Univ. Press, Cambridge (2018)
\bibitem{GS1}Galdi, G.P., and Silvestre A.L., Existence of time-periodic solutions to the Navier-Stokes equations around a moving body. {\em Pacific J. Math.}  {\bf 223}  251--267   (2006)
\bibitem{GaSo}Galdi, G.P., and Sohr, H., Existence and uniqueness of time-periodic physically reasonable Navier--Stokes flow past a body. {\em Arch. Ration. Mech. Anal.}  {\bf 172} 363--406   (2004)
\bibitem{GH}Geissert, M., Hieber, M., and Huy, N.-T., A general approach to time periodic incompressible viscous fluid flow problems. {\em Arch. Ration. Mech. Anal.}  {\bf 220} 1095--1118 (2016)
\bibitem{Hey} Heywood, J.G.,  
The Navier-Stokes equations: on the existence, regularity and decay of solutions, {\it Indiana Univ. Math. J.}, {\bf 29}, 639--681 (1980)
\bibitem{HG}Hieber, M., Huy, N.-T., and  Seyfert, A., On periodic and almost periodic solutions to incompressible viscous fluid flow problems on the whole line.  {\em Mathematics for nonlinear phenomena—analysis and computation}, 51-81, Springer Proc. Math. Stat., {\bf 215}, Springer, Cham (2017)
\bibitem{H}Hishida, T., Large time behavior of a generalized Oseen evolution operator, with applications to the Navier--Stokes flow past a rotating obstacle. {\em Math. Ann. }  {\bf 372}  915--949 (2018)
\bibitem{NTH}Huy, N.-T., Periodic motions of Stokes and Navier-Stokes flows around a rotating obstacle. {\em Arch. Ration. Mech. Anal.}  {\bf 213}  689--703 (2014)
\bibitem{KMT}Kang, K.,  Miura, H., and
Tsai, T--P.,
Asymptotics of small exterior Navier--Stokes flows with non-decaying boundary data. {\em 
Comm. Partial Differential Equations}  {\bf 37} 1717--1753  (2012)
\bibitem{KS}Korolev, A., and \v{S}ver\'ak, V.,  On the large-distance asymptotics of steady
state solutions of the Navier--Stokes equations in 3D exterior domains, {\em Ann.
I. H. Poincar\'e}, {\bf 28} 303--313 (2011)

\bibitem{KoNo}Kozono, H., and
Nakao, M., 
Periodic solutions of the Navier--Stokes equations in unbounded domains.
{\em Tohoku Math. J.} {\bf 48}  33--50 (1996)
\bibitem{Kyed}Kyed, M., The existence and regularity of time-periodic solutions to the three-dimensional Navier--Stokes equations in the whole space. {\em Nonlinearity}  {\bf 27}  2909--2935 (2014)
\bibitem{L}
Landau, L.
\newblock {A new exact solution of Navier-Stokes equations.}
\newblock {\em C. R. (Dokl.) Acad. Sci. URSS, n. Ser.}, {\bf 43} 286--288 (1944)
\bibitem{Ma} Maremonti, P.
Existence and stability of time-periodic solutions to the Navier-Stokes equations in the whole space. 
{\em Nonlinearity}  {\bf 4}  503--529 (1991)
\bibitem{Ri}Riley, N., Steady streaming, {\em Annu. Rev. Fluid Mech.} {\bf 33}, 43--65 (2001)
\bibitem{Sch}Schlichting, H., Boundary Layer Theory, 7th Edition, McGraw--Hill, New York (1979)
\bibitem{Sol}Solonnikov, V.A.,
     \newblock {Estimates of the solutions of the nonstationary Navier-Stokes system},
     \newblock {\em Zap. Naucn. Sem. Leningrad. Otdel. Mat. Inst. Steklov} (LOMI), \textbf{38} 153--201 (1973)

\bibitem{Y} Yamazaki, M., The Navier--Stokes equations in the weak--$L_n$   space with time-dependent external force. {\em Math. Ann.}  {\bf 317}  635--675 (2000)

\end{thebibliography}
